\documentclass{amsart}
\usepackage[french,english]{babel}
\usepackage{amsmath, amssymb}
\usepackage[mathscr]{euscript}
\usepackage[utf8]{inputenc}

\input xy
\xyoption{all}

\def\aut{\operatorname{Aut}}

\def\ab{\operatorname{ab}}
\def\Ab{\operatorname{Ab}}
\def\dcl{\operatorname{dcl}}
\def\dlog{\operatorname{dlog}}
\def\tp{\operatorname{tp}}
\def\C{\mathcal{C}}
\def\U{\mathcal{U}}

\def\acl{\operatorname{acl}}
\def\stp{\operatorname{stp}}
\def\dcf{\operatorname{DCF}_0}
\def\ccm{\operatorname{CCM}}

\def\Alb{\operatorname{Alb}}
\def\PP{\mathbb{P}}

\def\AA{\mathbb{A}}
\def\Gm{\mathbb{G}_{\operatorname{m}}}
\def\Ga{\mathbb{G}_{\operatorname{a}}}
\def\Aff{\operatorname{Aff}}
\def\PSL{\operatorname{PSL}}
\def\SL{\operatorname{SL}}

\def\Ind#1#2{#1\setbox0=\hbox{$#1x$}\kern\wd0\hbox to 0pt{\hss$#1\mid$\hss}
\lower.9\ht0\hbox to 0pt{\hss$#1\smile$\hss}\kern\wd0}
\def\ind{\mathop{\mathpalette\Ind{}}}
\def\Notind#1#2{#1\setbox0=\hbox{$#1x$}\kern\wd0\hbox to 0pt{\mathchardef
\nn=12854\hss$#1\nn$\kern1.4\wd0\hss}\hbox to
0pt{\hss$#1\mid$\hss}\lower.9\ht0 \hbox to
0pt{\hss$#1\smile$\hss}\kern\wd0}
\def\nind{\mathop{\mathpalette\Notind{}}}

\newtheorem{innercustomthm}{Theorem}
\newenvironment{customthm}[1]
  {\renewcommand\theinnercustomthm{#1}\innercustomthm}
  {\endinnercustomthm}

\newtheorem*{theorem*}{Theorem}

\newtheorem{theorem}{Theorem}[section]
\newtheorem{proposition}[theorem]{Proposition}
\newtheorem{lemma}[theorem]{Lemma}
\newtheorem{corollary}[theorem]{Corollary}

\theoremstyle{definition}
\newtheorem{definition}[theorem]{Definition}
\newtheorem{fact}[theorem]{Fact}

\theoremstyle{remark}
\newtheorem{remark}[theorem]{Remark}

\title[Abelian reductions]{Abelian reduction in differential-algebraic and bimeromorphic geometry}

\author{R\'emi Jaoui}
\address{R\'emi Jaoui\\
Université Claude Bernard Lyon 1, CNRS UMR 5208, Institut Camille Jordan, 43 Blvd. du 11 novembre 1918, 69622 Villeurbanne, France}
\email{jaoui@math.univ-lyon1.fr}

\author{Rahim Moosa}
\address{Rahim Moosa\\
University of Waterloo\\
Department of Pure Mathematics\\
200 University Avenue West\\
Waterloo, Ontario \  N2L 3G1\\
Canada}
\email{rmoosa@uwaterloo.ca}


\date{\today}

\begin{document}
\selectlanguage{english} 

\begin{abstract}
A new tool for the model theory of differentially
closed fields and of compact complex manifolds is here developed. In such settings, it is shown that a type internal to the field of constants (resp. to the
projective line) admits a maximal image whose binding group is an
abelian variety. The properties of such {\em abelian reductions} are
investigated in the Galois-theoretic framework provided by stability
theory.

Several geometric consequences for the birational geometry of algebraic
vector fields of characteristic zero are then deduced. In particular,
\begin{enumerate}
\item it is shown that if some cartesian power of an algebraic vector
field admits a nontrivial rational first integral then already the
second power does,
\item two-dimensional isotrivial algebraic vector fields are classified
up to \linebreak birational equivalence and
\item algebraic vector fields whose finite covers admit no nontrivial
factors are studied in arbitrary dimension.
\end{enumerate}
Analogues of these results in bimeromorphic geometry are also obtained.
\end{abstract}

 \maketitle
\selectlanguage{french} 
\begin{abstract}
\noindent
Un nouvel outil pour la théorie des modèles des corps différentielle-ment clos et des variétés compactes complexes est développé dans cet article. Dans ces théories, il est montré qu'un type interne au corps des constantes (resp. à la droite projective) admet une plus grande image dont le groupe de liaison est une variété abélienne. Les propriétés de ces {\em{reductions abéliennes}} sont étudiées dans le formalisme de la théorie géométrique de la stabilité. 

Plusieurs conséquences concernant la géométrie birationnelle des champs de vecteurs algébriques de caractéristique zéro sont alors obtenues. En particulier, 
\begin{itemize}
\item[(1)] il est montré que si une puissance cartésienne d'un champ de vecteurs algébrique admet une intégrale première rationnelle non triviale alors la seconde puissance cartésienne vérifie déjà cette propriété,

\item[(2)] les champs de vecteurs algébriques isotriviaux de dimension deux sont classifiés à équivalence birationnelle près, 

\item[(3)] les champs de vecteurs algébriques dont tous les recouvrement finis n'admettent aucun facteur non trivial sont étudiés en dimension arbitraire.
\end{itemize}
Des résultats analogues en géométrie biméromorphe sont aussi obtenus. 
\end{abstract}     

\maketitle

\selectlanguage{english} 

\vfill
\pagebreak

\setcounter{tocdepth}{1}
\tableofcontents

\section{Introduction}

\noindent
This paper is an example of the role model theory often plays of transferring ideas from one geometric context to another, of recognising and facilitating analogies between seemingly disparate areas of mathematics.
In this case, we are inspired by a notion in bimeromorphic geometry (Fujiki's {\em relative Albanese}), which we abstract and develop in the general setting of finite rank geometric stability theory, introducing thereby a new tool to the subject (namely {\em abelian reductions}).
We then apply the results we obtain back to bimeromorphic geometry, but also, significantly, to the birational geometry of algebraic vector fields.

Let us begin with our applications to differential-algebraic geometry, of which there are four.
When applicable, we state our results both geometrically (about algebraic vector fields) and model-theoretically (about the first order theory of differentially closed fields of characteristic zero -- $\dcf$).
Fix an algebraically closed field $k$ of characteristic zero and an irreducible algebraic vector field $(X,v)$ over $k$.
We denote the cartesian product of vector fields $(X,v)$ and $(Y,w)$ by $(X\times Y, v\times w)$, and cartesian powers by $(X^n,v^n)$.

\begin{theorem}
\label{nwodeg-geometric-intro}
Suppose $(X^2,v^2)$ admits no nontrivial rational first integral.
Then neither does $(X^n,v^n)$, for all $n\geq 1$.
In fact, if $(Y,w)$ is any other irreducible algebraic vector field over $k$ without nontrivial rational first integrals, then the product $(X\times Y,v\times w)$ admits no nontrivial rational first integral.
\end{theorem}

\begin{customthm}{\ref{nwodeg-geometric-intro} (Model-theoretic formulation)}
\label{nwodeg-intro}
Suppose $p$ is a complete stationary type of finite rank over constant parameters in $\dcf$.
If $p$ is not orthogonal to the constants then $p^{(2)}$ is not weakly orthogonal to the constants.
\end{customthm}

Recall that an algebraic vector field having no nontrivial rational first integrals is equivalent to the induced derivation on the function field having no new constants.
See Section~\ref{sec:dcf} below for details, where this theorem appears as Theorem~\ref{nwodeg} and Corollary~\ref{nwodeg-geometric}.
Theorem~\ref{nwodeg-geometric-intro} is of course sharp, in the sense that there are many algebraic vector fields which admit no nontrivial rational first integrals but whose cartesian squares do; the vector field on $\AA^1$ that is identically equal to $1$ is such an example.
Regarding the model-theoretic formulation, it is well known that if~$p$ is nonorthogonal to the constants then some Morley power $p^{(n)}$ is not weakly orthogonal to the constants; our result says that one can take $n=2$ in the {\em autonomous case}, that is, when the parameters live in the field of constants.
This assumption that $p$ be autonomous is necessary; in general there can be no bound independent of the dimension of $p$, see the discussion following Theorem~\ref{nwodeg} below.

Theorem~\ref{nwodeg-geometric-intro} can also be interpreted as an algebraic analogue of  a well known result in the ergodic theory of flows and discrete dynamical systems.
If $X$ is a measurable space and $\phi: X \rightarrow X$ is a measurable map (or a measurable flow) preserving a probability measure $\mu$ then the dynamical system $(X,\phi,\mu)$ is {\em ergodic} if every measurable first integral of $(X,\phi)$ is constant $\mu$-almost everywhere (or equivalently if the only elements $f$ of the Hilbert space $L^2(X,\mu)$ satisfying $f \circ \phi = f$ are the constants).
Similarly to Theorem~\ref{nwodeg-geometric-intro}, one has that if $(X,\phi,\mu) \times (X,\phi,\mu)$ is ergodic then $(X,\phi,\mu) \times  (Y,\psi,\nu)$ is ergodic for every ergodic system $(Y,\psi,\nu)$.
See, for example, \cite[Proposition~4.4]{furstenberg}.
The model-theoretic formulation of Theorem~\ref{nwodeg-geometric-intro} can, in turn, be seen as an algebraic strengthening of the dynamical criterion for orthogonality to the constants given by the first author in~\cite{jaoui20}. 

Our second application is an explicit birational classification of algebraic vector fields of dimension~2 that are {\em isotrivial} in the sense that they are isomorphic to a trivial vector field, though possibly after some base extension by a differential field containing $k$.
Again, see Section~\ref{sec:dcf} for details, where this appears as Theorem~\ref{dim2} and Corollary~\ref{classification-ordertwo}.

\begin{theorem}
\label{dim2-intro}
Suppose $\dim X=2$, $(X,v)$ is isotrivial, and $(X,v)$ admits no nontrivial rational first integrals.
Then one of the following holds:
\begin{enumerate}
 \item
 $(X,v)$ is birationally equivalent over $k$ to
 \begin{itemize}
 \item[(a)]
 the vector field on $\mathbb A^2$ given by $w(x,y)=(x,y,\lambda x, \mu y)$ for some $\mathbb Z$-linearly independent  $\lambda,\mu\in k$, or
 \item[(b)]
 the vector field on $\mathbb A^2$ given by $w(x,y)=(x,y,\lambda, \mu y)$ for some nonzero $\lambda,\mu\in k$, or,
\item[(c)]
a translation invariant vector field on an abelian surface, or
\end{itemize}
\item
$X$ admits a dominant rational map $X\to E$ to an elliptic curve $E$ over $k$ such that $v$ lifts a translation invariant vector field $w$ on $E$, and such that after localisation to an \'etale open set $U\to E$, $(X,v)$ becomes birationally equivalent to 
the projectivisation of a vector field on $\AA^2 \times U$ that lifts $w_U$ and is linear on the fibres over $U$.
\end{enumerate}
\end{theorem}

\begin{customthm}{\ref{dim2-intro} (Model-theoretic formulation)}
Suppose, in $\dcf$, that $p$ is a complete stationary type of dimension~$2$ over constant parameters that is internal and weakly orthogonal to the constants.
Then one of the following two cases occurs:
\begin{enumerate}
\item
$p$ is interdefinable with the generic type of a logarithmic-differential equation on either $\mathbb G_m^2$, $\mathbb G_m\times\mathbb G_a$, or an abelian surface; or
\item
there is a fibration $f:p\to q$ where $q$ is the generic type of a logarithmic-differential equation on an elliptic curve whose fibre is interdefinable with the generic type of a Riccati equation.
\end{enumerate}
\end{customthm}

In the spirit of analogies, Theorem~\ref{dim2-intro} should be seen as a differential-algebraic counterpart to Fujiki's~\cite{fujiki83} classification of $3$-dimensional compact K\"ahler manifolds of algebraic dimension~$1$.

Our third application, appearing as Theorem~\ref{nofibre-dcf} and Corollary~\ref{nofibre-dcf-geometric} below, provides a structural dichotomy for algebraic vector fields all of whose finite covers admit no nontrivial factors:

\begin{theorem}
\label{nofibre-dcf-geometric-intro}
Assume that, for every generically finite cover $(Y,w)\to (X,v)$, $(Y,w)$ admits no dominant morphism $(Y,w)\to (Z,u)$ with $0<\dim Z<\dim Y$.
Then either:
\begin{itemize}
\item[(i)]
for all algebraic vector fields $(Y,w)$ over $k$, every proper invariant subvariety of $(X\times Y, v\times w)$ over $k$ that projects dominantly onto both $X$ and $Y$ is a generically finite cover of $(Y,w)$, or
\item[(ii)]
$(X,v)$ is a generically finite cover of a translation invariant vector field on a simple abelian variety over $k$ of dimension $>1$.
\end{itemize}
\end{theorem}

\begin{customthm}{\ref{nofibre-dcf-geometric-intro} (Model-theoretic formulation)}
Suppose $p=\tp(a/B)$ is a complete stationary nonalgebraic type over constant parameters that satisfies the following form of exchange: if $c\in\acl(Ba)\setminus\acl(B)$ then $a\in\acl(Bc)$.
Then $p$ is either minimal or admits a finite-to-one definable function $f:p\to q$ where $q$ is the generic type of a logarithmic-differential equation on a simple abelian variety of dimension $>1$.
\end{customthm}

Again we have a guiding precedent to Theorem~\ref{nofibre-dcf-geometric-intro} in bimeromorphic geometry: the work of Campana, Oguiso, and Peternell on ``minimal fibrations" of compact K\"ahler manifolds, see~\cite[Theorem~2.4]{COP} as well as its model-theoretic account in~\cite{moosa2014some}.
The exchange principle appearing in the model-theoretic formulation of Theorem~\ref{nofibre-dcf-geometric-intro} strengthens the condition of having {\em no proper fibrations}, which was introduced and studied in~\cite{moosa2014some}, by replacing $\dcl$ with $\acl$.
Minimal types certainly satisfy this strengthened condition, and 
Theorem~\ref{nofibre-dcf-geometric-intro} describes explicitly all the nonminimal autonomous examples in $\dcf$.
The result fails if we move outside the autonomous setting.

Finally, in a somewhat different direction, we obtain a strong constraint on which projective algebraic varieties can admit an algebraic vector field having no nontrivial rational first integrals.
(Note that by a theorem of Buium~\cite{Buium}, any algebraic vector field on a projective variety is necessarily isotrivial.)
The following theorem, which appears as Theorem~\ref{albanese-surjective} below, can be viewed as a differential-algebraic analogue of a special case of a theorem of Fujiki~\cite[Theorem~1]{fujiki83} on compound Moishezon fibrations of compact K\"ahler manifolds.

\begin{theorem}
\label{albanese-surjective-intro}
Suppose $Y$ is a smooth  projective variety over $k$ that admits an algebraic vector field with no nontrivial rational first integrals.
Then the Albanese morphism $\alpha:Y\to\Alb(Y)$ is surjective with generic fibre absolutely irreducible and birationally equivalent to an algebraic homogeneous space for a connected linear algebraic group.
\end{theorem}

We now describe the fundamental abstract result that underlies all these applications.
We use here freely notions from geometric stability theory, especially around internality, binding groups, and orthogonality, but that material is reviewed in Section~\ref{sec:prelims} below.
Recall that the Chevalley decomposition describes an algebraic group as the extension of an abelian variety by a linear algebraic group.
The following theorem, which combines Theorems~\ref{abredexists} and~\ref{abred=alb} below, describes how the Chevalley decomposition of the binding group of an internal type is reflected in the structure of the type itself.

\begin{theorem}
\label{abred-intro}
Suppose $\U$ is a sufficiently saturated model of a complete totally transcendental theory admitting elimination of imaginaries, with a $0$-definable pure algebraically closed field $\C$.
Suppose $p\in S(B)$ is a complete stationary type over a set of parameters $B$.
If $p$ is $\C$-internal then there exists $p_{\ab}\in S(B)$ and a fibration $f:p\to p_{\ab}$ such that:
\begin{enumerate}
\item
the binding group of $p_{\ab}$ relative to $\C$ is definably isomorphic to the $\C$-points of an abelian variety over $\C$, and
\item
if $g:p\to q$ is any definable function with $q\in S(B)$ having binding group relative to~$\C$ definably isomorphic to the $\C$-points of an abelian variety over~$\C$ then $g$ factors through $f$.
\end{enumerate}
We call $f:p\to p_{\ab}$ the {\em abelian reduction} of $p$.
If, in addition, $B=\acl(B)$ and $p$ is weakly $\C$-orthogonal, then the binding group of $p_{\ab}$ is the abelian part in the Chevalley decomposition of the binding group of $p$.
\end{theorem}

Our applications to algebraic vector fields (namely Theorems~\ref{nwodeg-geometric-intro} through~\ref{albanese-surjective-intro}) are deduced by investigating the consequences of Theorem~\ref{abred-intro} specialised to the theory of differentially closed fields, with $\C$ the field of constants, combined with the failure of the inverse differential Galois problem over constant parameters; namely the fact that the only linear algebraic groups that can arise as binding groups of autonomous types are the commutative ones (see Fact~\ref{linearbg} below).

Theorem~\ref{abred-intro} can also be specialised to the theory of compact complex manifolds ($\ccm$) with $\C$ (an elementary extension of the) complex field living on the projective line.
Moreover, the analogue of the failure of the inverse Galois problem over constant parameters holds also for compact K\"ahler manifolds; in fact, no nontrivial linear algebraic group can arise as the binding group of an autonomous type.
This is Proposition~\ref{linearbg-ccm} below, and is deduced from a theorem of Fujiki on the relative Albanese.
Applications to compact K\"ahler manifolds, analogous to Theorems~\ref{nwodeg-geometric-intro} through~\ref{albanese-surjective-intro} above, follow.
However, many of the results one obtains in this way were already known; indeed, they are the theorems we have mentioned as inspiring our work.
But the following consequences are, to the best of our knowledge, new results in bimeromorphic geometry:

\begin{theorem}
\label{nwodeg-ccm-intro}
Suppose $X$ is a compact complex variety of K\"ahler-type, $Y$ is Moishezon, and $f:X\to Y$ is a fibre space.
\begin{itemize}
\item[(a)]
If $f^{(2)}:X_Y^{(2)}\to Y$ induces an isomorphism of meromorphic function fields then so does $f^{(n)}:X_Y^{(n)}\to Y$ for all~$n\geq 1$.
\item[(b)]
If every irreducible proper analytic subset of $X_Y^{(2)}$ that projects onto $X$ in both co-ordinates is of dimension $\dim X$, then for any surjective morphism $g:Z\to Y$, any irreducible proper analytic subset $W\subseteq X_{(Z)}$ that projects onto both $X$ and $Z$ must be of dimension $\dim Z$.
\end{itemize}
\end{theorem}

\begin{customthm}{\ref{nwodeg-ccm-intro} (Model-theoretic formulation)}
Suppose that $p$ is a complete stationary type in $\ccm$ over constant parameters living on the sort of a compact complex variety of K\"ahler-type.
\begin{itemize}
\item[(a)]
If $p$ is not orthogonal to the constants then $p^{(2)}$ is not weakly orthogonal to the constants.
\item[(b)]
The degree of nonminimality of $p$ is at most~$1$.
In particular, if every pair of distinct realisations of $p$ is independent then $p$ is minimal.
\end{itemize}
\end{customthm}

This is done in Section~\ref{sec:ccm} below, in particular, in Corollary~\ref{autccm}, Theorem~\ref{autccm-geometric}, and Theorem~\ref{nwodeg-ccm}.
We point the reader to the beginning of that section for a proper explanation of the terminology being used,
but let us give some brief explanations now: by $f^{(n)}:X_Y^{(n)}\to Y$ we mean the fibration whose general fibres are the $n$th cartesian powers of the fibres of $f$, and by $X_{(Z)}$ we mean the unique maximal dimensional irreducible component of the fibred product $X\times_YZ$ that projects onto~$Z$.
The assumption in the above theorem that $Y$ is Moishezon, namely bimeromorphic to a complex projective algebraic variety, reflects in this setting the fact that we are concerned with autonomous types.

Theorem~\ref{nwodeg-ccm-intro}(a) is the analogue of Theorem~\ref{nwodeg-geometric-intro} for compact complex manifolds.

Regarding the model-theoretic formulation of Theorem~\ref{nwodeg-ccm-intro}(b), the {\em degree of nonminimality} was introduced by the second author and James Freitag in~\cite{nmdeg}, where it was shown to be bounded by the $U$-rank of $p$ plus~$1$.
This was improved in~\cite{nmdeg<2} where an absolute upper bound of~$2$ was obtained, both in $\dcf$ and $\ccm$.
Moroever, in the case of $\dcf$ it was shown that autonomous types have degree of nonminimality at most~$1$.
Theorem~\ref{nwodeg-ccm-intro}(b) is an extension of that result to~$\ccm$.

\medskip
\subsection{Acknowledgements}
The first author was partially supported by the ANR-DFG program GeoMod (Project number 2100310201).
The second author was partially supported by an NSERC DG.
The second author would also like to thank the hospitality of the Fields Institute where some of this work was carried out during the 2021 Thematic Programme on Trends in Pure and Applied Model Theory.

\medskip
\subsection{Notation}
When dealing with sets of parameters $B,C$, in a first-order structure, we will employ the common abuse of notation in model theory where $BC$ stands for $B\cup C$; and if $a=(a_1,\dots,a_n)$ then $Ba$ stands for $B\cup\{a_1,\dots,a_n\}$.

\bigskip
\section{Preliminaries on binding groups}
\label{sec:prelims}

\noindent
We work throughout this paper, unless stated otherwise, in a fixed sufficiently saturated model $\U$ of a complete totally transcendental theory admitting elimination of imaginaries, with a $0$-definable pure algebraically closed field $\C$.
Also unless stated otherwise, ``definable" always means ``definable with parameters".

Given any parameter set $B$ we set
$$\aut_B(\U/\C):=\{\sigma\in\aut(\U): \sigma \text{ fixes $B\C$ pointwise}\},$$
and given $p\in S(B)$, we set
$$\aut(p/\C):=\{\sigma\upharpoonright_{p(\U)}:\sigma\in\aut_B(\U/\C)\}.$$
That is, $\aut(p/\C)$ is the group of elementary permutations of $p(\U)$ that fix $B\C$ pointwise.

Recall that we say $p$ is {\em $\C$-internal} if it is stationary and there exists $B'\supseteq B$ such that $p(\U)\subseteq\dcl(B'\C)$.
If we replace $\dcl$ by $\acl$ in the above definition then we get the notion of {\em almost $\C$-internal}.
The binding group theorem tells us that if $p$ is $\C$-internal then $\aut(p/\C)$, along with its action on $p(\U)$, is isomorphic to a $B$-definable group acting relatively $B$-definably on $p(\U)$.
Note that, by definition, this is always a faithful group action.
Unless stated otherwise, we identify $\aut(p/\C)$ with this $B$-definable group and call $\aut(p/\C)$ the {\em binding group} of $p$.

A stationary type $p\in S(B)$ is {\em weakly $\C$-orthogonal} if for every (equivalently some) $a\models p$, and every finite tuple $c$ from $\C$, we have that $a\ind_Bc$.
If $p$ is both $\C$-internal and weakly $\C$-orthogonal then the action of $\aut(p/\C)$ on $p(\U)$ is transitive and the type $p$ itself is isolated.
Hence, in this case, we have that $(\aut(p/\C),p(\U))$ is a $B$-definable {\em homogeneous space}; that is, a faithful and transitive group action.
In fact, if $p$ is $\C$-internal then weak $\C$-orthogonality,  transitivity of the binding group action, and isolation of $p$, are all equivalent.

Every $\C$-internal type $p\in S(B)$ admits a {\em fundamental system of realisations}, namely a finite tuple of realisations $\eta$ of $p$ such that $p(\U)\subseteq\dcl(B\eta\C)$.
That is, in the definition of internality we can take $B'=B\eta$.
The action of $\aut(p/\C)$ on $p(\U)$ is then determined by its action on $\eta$.
A useful fact is that we can always find a fundamental system of realisations satisfying $p^{(n)}$ for some $n\geq 1$.
Here $p^{(n)}$ denotes the Morley power of $p$, namely the type of an independent sequence of realisations.

We will sometimes work over algebraically closed sets of parameters.
Let us record right away how the binding group behaves with respect to passing to the algebraic closure.

\begin{lemma}
\label{autbar}
Suppose $p\in S(B)$ is $\C$-internal  and let $\overline p\in S(\acl(B))$ be the unique extension of $p$.
Then $\aut(\overline p/\C)$ is a finite index normal definable subgroup of $\aut(p/\C)$.
\end{lemma}

\begin{proof}
As abstract permutation groups of $p(\U)=\overline p(\U)$, note that $\aut(\overline p/\C)$ is a subgroup of $\aut(p/\C)$; it is the subgroup of permutations that have a lifting to $\aut_{\acl(B)}(\U/\C)$.
But we wish to identify $\aut(\overline p/\C)$ definably as a subgroup of $\aut(p/\C)$.
This can be done canonically by going back to the definable structure given to the binding groups, but here is a (non-canonical) short cut:
Fix a fundamental system of realisations $\eta$ for $p$.
Then it is also a fundamental system for $\overline p$.
Now consider the function $\iota:\aut(\overline p/\C)\to\aut(p/\C)$ given by $\iota(\sigma)=\tau$ if the action of $\sigma$ and $\tau$ agree on $\eta$.
Since both $\aut(\overline p/\C)$ and $\aut(p/\C)$ act relatively definably, and faithfully, this is a definable (over $B\eta$)  injective homomorphism.

Since $\aut_{\acl(B)}(\U/\C)$ is a normal subgroup of $\aut_B(\U/\C)$, it follows that $\aut(\overline p/\C)$ is a normal subgroup of $\aut(p/\C)$.
Moreover, the restriction of $\aut_B(\U/\C)\to\aut(p/\C)$ induces a surjective homomorphism
$$\aut_B(\U/\C)/\aut_{\acl(B)}(\U/\C)\to \aut(p/\C)/\aut(\overline p/\C).$$
But
$\aut_B(\U/\C)/\aut_{\acl(B)}(\U/\C)$
is isomorphic to the group of permutations of $\acl(B)$ that lift to automorphisms in $\aut_B(\U/\C)$, and is therefore profinite.
So $\aut(p/\C)/\aut(\overline p/\C)$ is a profinite group.
As it is interpretable in a totally transcendental theory, it must be finite.
\end{proof}

The binding group itself can be definably embedded in a cartesian power of~$\C$.
Because the induced structure on~$\C$ is that of a stably embedded pure algebraically closed field, the binding group is definably isomorphic with the $\C$-points of an algebraic group.
In general, this definable isomorphism may require additional parameters.
Let us record for future use the well known fact that when the binding group is commutative then no additional parameters are required:

\begin{lemma}
\label{overB}
Suppose $p \in S(B)$ is $\mathcal C$-internal and that $\aut(p/\C)$ is commutative.
Then there is an algebraic group $H$ over $\dcl(B)\cap\C$ and a $B$-definable isomorphism between $\aut(p/\C)$ and $H(\C)$.
\end{lemma}
 
\begin{proof} 
Set $G:=\aut(p/\C)$ and denote by $\sigma \mapsto g_\sigma$ the map $\aut_B(\U/\mathcal C) \rightarrow G$ given by restriction to $p(\U)$.

We first show that $G \subset\dcl(B\mathcal C)$.
It is enough to see that $\sigma(g) = g$ for all $g \in G$, $\sigma \in \aut_B(\U/ \mathcal C)$.
Consider $a \models p$.
Then $ga$ and $a$ realise $p$ so $\sigma(ga) = g_\sigma(ga)$ and $\sigma(a) = g_\sigma a$.
It follows that we can compute $\sigma(ga)$ in two different ways: 
$$ \sigma(ga) = \sigma(g)\sigma(a) = \sigma(g)(g_\sigma a) = (\sigma(g)g_\sigma)a $$
and 
$$ \sigma(ga) =  g_\sigma(ga) = (g_\sigma g)a.$$
It follows that the equality $(g_\sigma g)a = (\sigma(g)g_\sigma)a $ holds for every realisation $a$ of $p$. Since the action of the binding group $G$ on $p(\U)$ is faithful, we obtain that  $\sigma(g)g_\sigma = g_\sigma g$,
and hence by commutativity of $G$ that $\sigma(g)=g$ as desired.

By compactness it follows that $G$ is the image of a $B$-definable subset of some cartesian power of $\C$ under a $B$-definable function.
Since $\C$ is a stably embedded pure algebraically closed field, it eliminates imaginaries, and we have a $B$-definable isomorphism $f:G\to D$ where $D$ is a $(\dcl(B)\cap\C)$-definable group in $(\C,0,1,+,-,\times)$.
But every such group is of the form $H(\C)$ for some algebraic group $H$ over $\dcl(B)\cap\C$.
\end{proof}

\medskip
\subsection{Galois correspondence}
A $\C$-internal type $p\in S(B)$ is {\em fundamental} if any (equivalently some) single realisation of $p$ is a fundamental system.
Note that if $p$ is fundamental $\C$-internal and weakly $\C$-orthogonal, then $(\aut(p/\C),p(\U))$ is a principal homogeneous space.
In any case, in such a situation we have the following Galois correspondence, see~\cite[Appendix~B]{udicomputing}.

\begin{fact}
\label{gcorr}
Suppose $p=\tp(a/B)$ is fundamental $\C$-internal and weakly $\C$-orthogonal.
For any set $E\subseteq\dcl(Ba)$ let
$$H_E:=\{\sigma\in\aut(p/\C):\text{ some (any) lifting of $\sigma$ to $\aut_B(\U/\C)$ fixes $E$ pointwise}\}.$$
\begin{itemize}
\item[(a)]
For any $\dcl(B)\subseteq E\subseteq\dcl(Ba)$ definably closed, $H_E$ is a $Ba$-definable subgroup of $\aut(p/\C)$.
Every $Ba$-definable subgroup arises in this way.
\item[(b)]
If $E=\dcl(Bb)\subseteq\dcl(Ba)$ and $q:=\tp(b/B)$ is fundamental then $H_E$ is a $B$-definable normal subgroup of $\aut(p/\C)$.
Every $B$-definable subgroup of $\aut(p/\C)$ is normal and arises in this way.
Moreover, in this case $\aut(q/\C)$ is $B$-definably isomorphic to the quotient $\aut(p/\C)/H_E$ along with its natural action on $q(\U)$.
\end{itemize}
\end{fact}

While the assumption of fundamentality in the Galois correspondence is essential, we can always obtain a fundamental type from a $\C$-internal one:
Fixing $\eta$ a  fundamental system of realisations for $p\in S(B)$ and assuming $B=\acl(B)$ to ensure stationarity, we obtain that $q:=\tp(\eta/B)$ is fundamental $\C$-internal.
Moreover, $\aut(p/\C)=\aut(q/\C)$.
What we mean by this, of course, is that the map $\sigma\mapsto\widehat\sigma\upharpoonright_{q(\U)}$ is an isomorphism from $\aut(p/\C)$ to $\aut(q/\C)$, where $\widehat\sigma\in\aut_B(\U/\C)$ is a lifting of $\sigma$.
That this is well-defined follows from the fact that $q(\U)\subseteq p(\U)^n$.
That it is injective follows from the fact that $p(\U)\subseteq\dcl(B\eta\C)$.
Surjectivity is clear: given $\tau\in\aut(q/\C)$ lift it to some $\widehat \tau\in\aut_B(\U/\C)$ and observe that $\sigma:=\widehat\tau\upharpoonright_{p(\U)}\in\aut(p/\C)$ maps to $\tau$.

If, in addition, we know that $p$ is weakly $\C$-orthogonal, then we can choose our fundamental system $\eta$ so that $q=\tp(\eta/B)$ is also weakly $\C$-orthogonal.
Indeed, in this case $p$ is isolated and the fact that $p(\U)\subseteq\dcl(B\eta\C)$ can be expressed as a first order property of~$\eta$.
It follows that we can choose such an $\eta$ in a prime model over~$B$.
Hence $q$ is isolated, and this forces weak $\C$-orthogonality.

Let us record the above observations as a fact for future use:

\begin{fact}
\label{makefundamental}
Suppose $B=\acl(B)$ and $p\in S(B)$ is $\C$-internal and weakly $\C$-orthogonal.
Then there exists a fundamental system of realisations, $\eta$ of $p$, such that $q:=\tp(\eta/B)$ is fundamental $\C$-internal, weakly $\C$-orthogonal, and $\aut(p/\C)=\aut(q/\C)$.
\end{fact}

Here is one useful consequence:

\begin{corollary}
\label{connectedbind}
Suppose $B=\acl(B)$ and $p\in S(B)$ is $\C$-internal and weakly $\C$-orthogonal.
Then $\aut(p/\C)$ is connected.
\end{corollary}

\begin{proof}
Let $G:=\aut(p/\C)$.
Let $\eta$ be a fundamental system of realisations given by Fact~\ref{makefundamental}, and set $q=\tp(\eta/B)$.
Then $G=\aut(q/\C)$ also.
Let $G^\circ$ be the connected component of~$G$.
This is a normal $B$-definable subgroup, and hence, by the Galois correspondence applied to $q$, namely Fact~\ref{gcorr}(b), there exists $b\in\dcl(B\eta)$, such that $G/G^\circ$ is the binding group of $r:=\tp(b/B)$.
But $r$ is also weakly $\C$-orthogonal, and hence the action of $G/G^\circ$ is transitive on $r(\U)$.
This forces $r(\U)$ to be finite, and as $B=\acl(B)$ we must have that $r$ is realised in $B$.
But then $\aut(r/\C)$ is trivial.
Hence $G^\circ=G$, as desired.
\end{proof}

\medskip
\subsection{Fibrations}
Given $p, q\in S(B)$ stationary, by a {\em definable function $f:p\to q$} we mean that $f$ is a $B$-definable function and $f(a)\models q$ for some (equivalently all) $a\models p$.
If $a\in\acl(Bf(a))$ then we say that $f$ is {\em finite-to-one}.
We say that $p$ and $q$ are {\em interdefinable} if there are definable functions $p\to q$ and $q\to p$.
By a {\em fibration} we mean a definable function $f:p\to q$ such that $\tp(a/Bf(a))$ is stationary for some (equivalently all) $a\models p$.
By a {\em fibre} of $f:p\to q$ we mean a type of the form $\tp(a/Bf(a))$.
Note that all the fibres are $B$-conjugate.

Taking Morley products preserves the property of being a fibration:

\begin{lemma}
\label{Morleyfibration}
Suppose $f_i:p_i\to q_i$ are fibrations for $i=1,2$.
Then the coordinate-wise product $f_1\times f_2:p_1\otimes p_2\to q_1\otimes q_2$ is a fibration.
\end{lemma}

\begin{proof}
Clearly,  $f_1\times f_2:p_1\otimes p_2\to q_1\otimes q_2$ is a definable function.
We need only show that if $\eta:=(a_1,a_2)\models p_1\otimes p_2$ then $\tp(\eta/Bf(\eta))$ is stationary.
This is straightforward forking calculus.
Let $b_i:=f_i(a_i)$ for $i=1,2$ and $C:=\acl(Bb_1b_2)$.
We need to show that $\tp(a_1a_2/Bb_1b_2)$ determines $\tp(a_1a_2/C)$.

Let $(a_1',a_2')\models\tp(a_1,a_2/Bb_1b_2)$.
Note that $f_1(a_1')=b_1$ and $f_2(a_2')=b_2$.
Now, from $a_1\ind_Bb_2$ and $a_1\ind_{Bb_1b_2}C$ we get that $a_1\ind_{Bb_1}C$.
Similarly, $a_1'\ind_{Bb_1}C$.
By stationarity of $\tp(a_1/Bb_1)$, which is the fact that $f_1:p_1\to q_1$ is a fibration, we get that $\tp(a_1/C)=\tp(a_1'/C)$.
Replacing everything by a $C$-conjugate, we may as well assume that $a_1'=a_1$.
That is, $\tp(a_2'/Ba_1b_2)=\tp(a_2/Ba_1b_2)$.
It follows that both $a_2$ and $a_2'$ are independent of $a_1$ over $Bb_2$.
Hence, by algebraicity, both are independent of $Ca_1$ over $Bb_2$.
Stationarity of $\tp(a_2/Bb_2)$ therefore implies that $\tp(a_2/Ca_1)=\tp(a_2'/Ca_1)$.
Hence, $(a_1',a_2')\models\tp(a_1,a_2/C)$, as desired.
\end{proof}

Assuming $p$ is $\C$-internal, a definable function $f:p\to q$ gives rise to a surjective $B$-definable homomorphism $f^\natural:\aut(p/\C)\to\aut(q/\C)$ of binding groups which can be described as follows: $f^\natural(\sigma)=\widehat\sigma\upharpoonright_{q(\U)}$ where $\widehat\sigma\in\aut(\U/\C)$ is any lifting of~$\sigma$.
If $f$ is finite-to-one then $\ker(f^\natural)$ will be finite.
See, for example,~\cite[Lemma~3.1]{jin-moosa}.
When $f$ is a fibration we can describe $\ker(f^\natural)$ as follows.
(See also~\cite[Corollary~5.12]{leo}.)

\begin{lemma}
\label{kernelstructure}
Suppose $p\in S(B)$ is a $\C$-internal type and $f:p\to q$ is a fibration.
Fix $n\geq 1$ and $\eta\models p^{(n)}$ such that $\eta$ is a fundamental system of realisations for $p$ and $f(\eta)$ is fundamental system of realisations for $q$.
Let $r:=\tp(\eta/Bf(\eta))$.
The kernel of $f^\natural:\aut(p/\C)\to\aut(q/\C)$ is definably isomorphic to $\aut(r/\C)$.
\end{lemma}

\begin{proof}
Note that $r$ is a fibre of $f:p^{(n)}\to q^{(n)}$, and hence, by Lemma~\ref{Morleyfibration}, $r$ is stationary.
Moreover, it too is $\C$-internal.
So $\aut(r/\C)$ is a $Bf(\eta)$-definable group acting (relatively) $Bf(\eta)$-definably on $r(\U)$.

Define $\rho:\aut(r/\C)\to\aut(p/\C)$ by $\sigma\mapsto\widehat\sigma\upharpoonright_{p(\U)}$ where $\widehat\sigma\in\aut_{Bf(\eta)}(\U/\C)$ is any lifting of $\sigma$.
This is well-defined because $\eta$ is a fundamental system for $p$ and hence any element of $\aut_B(\U/\C)$ that fixes $\eta$ is the identity on $p(\U)$.
Also, $\rho$ is definable (over $B\eta$) since $\rho(\sigma)=\tau$ if and only if $\sigma(\eta)=\tau(\eta)$.
As $r(\U)\subseteq p(\U)$ we have that $\rho$ is an embedding of $\aut(r/\C)$ in $\aut(p/\C)$.

We show that the image of $\rho$ is precisely $\ker(f^\natural)$.
Indeed, if $\sigma\in\aut(r/\C)$ then, as $\sigma$ lifts to an element of $\aut_{Bf(\eta)}(\U/\C)$, so does $\rho(\sigma)$.
Hence $f^\natural\rho(\sigma)$ fixes~$f(\eta)$.
But $f(\eta)$ is fundamental for $q$, and so this forces $f^\natural\rho(\sigma)=1$.
Conversely, if $\tau$ is in $\ker(f^\natural)$ then any lifting $\widehat\tau\in\aut_B(\U/\C)$ is contained in $\aut_{Bf(\eta)}(\U/\C)$ and hence $\sigma:=\widehat\tau\upharpoonright_{r(\U)}$ is contained in $\aut(r/\C)$ and $\rho(\sigma)=\tau$.
So $\rho$ is an isomorphism between $\aut(r/\C)$ and~$\ker(f^\natural)$.
\end{proof}

\medskip
\subsection{Definable Chevalley and Rosenlicht decompositions}
We record here the following definable versions of some important decomposition theorems for algebraic groups.
These definable versions are more or less well known and follow readily from the classical case using only the fact that $\C$ is a stably embedded pure algebraically closed field.
But we give some details for the sake of completeness.

\begin{fact}
\label{definablechevalley}
Suppose $G$ is a $B$-definable group that is definably isomorphic to the $\C$-points of a connected algebraic group over $\C$.
\begin{itemize}
\item[(a)]
There exists a unique maximal definable subgroup $L\leq G$ such that $L$ is definably isomorphic to the $\C$-points of a linear algebraic group over~$\C$.
Moreover, $L$ is normal, connected, $B$-definable, and $G/L$ is definably isomorphic to the $\C$-points of a connected abelian variety over $\C$.
\item[(b)]
There exists a unique minimal normal definable subgroup $D\leq G$ such that $G/D$ is definably isomorphic to the $\C$-points of a linear algebraic group over~$\C$.
Moreover, $G=D\cdot L$ and $D$ is also the unique minimal definable subgroup of $G$ with this property.
Moreover, $D$ is connected, commutative, and $B$-definable.
\end{itemize}
We call $L$ the {\em linear part}, $A:=G/L$ the {\em abelian part}, and $D$ the {\em anti-linear} part.
\end{fact}

\begin{proof}
Let $f:G\to \widehat G(\C)$ be a definable isomorphism where $\widehat G$ is a connected algebraic group over $\C$.
Note that we do not assume $f$ is defined over $B$.

The Chevalley decomposition of $\widehat G$ is a short exact sequence of connected algebraic groups over $\C$,
$$1\to\widehat L\to\widehat G\to\widehat A\to 0$$
where $\widehat L$ is linear and $\widehat A$ is abelian.
Moreover, $\widehat L$ is the unique maximal linear algebraic subgroup of $\widehat G$.
If we let $L:= f^{-1}(\widehat L(\C))$, then this is a definable subgroup of $G$ that is definably isomorphic to the $\C$-points of a linear algebraic group over~$\C$ (namely~$\widehat L$).
We show, first of all, that $L$ is maximal with this property, and uniquely so.
Suppose $N$ is a definable subgroup of $G$ that is definably isomorphic to the $\C$-points of a linear algebraic group over $\C$.
Then $f(N)$ is a definable subgroup of $\widehat G(\C)$.
As $\C$ is a stably embedded pure algebraically closed field, $f(N)$ must be of the form $\widehat N(\C)$ for some algebraic subgroup $\widehat N$ of $\widehat G$ over $\C$.
Moreover, $\widehat N$ must be linear since by assumption on $N$ we have that $\widehat N(\C)$ is definably isomorphic to the $\C$-points of some linear algebraic group over $\C$.
It follows that $\widehat N\leq\widehat L$, and hence
$$N=f^{-1}(\widehat N(\C))\leq f^{-1}(\widehat L(\C))=L,$$
as desired.
So $L$ is the unique maximal definably subgroup of $G$ that is definably isomorphic to the $\C$-points of a linear algebraic group over $\C$.
It follows that $L$ is in fact $B$-definable: any $B$-conjugate of $L$ would also be a maximal definable subgroup of $G$ that is definably isomorphic to the $\C$-points of a linear algebraic group over $\C$, and hence would be equal to $L$ by uniquenesss.
Finally, note that $G/L$ is definably isomorphic to $\widehat A(\C)$.
This proves part~(a).

The Rosenlicht decomposition of $\widehat G$ gives us a unique minimal normal algebraic subgroup $\widehat D$ such that the quotient is linear.
See~\cite[$\S$5]{brion} for details.
In particular, it is a fact that $\widehat D$ is also the unique minimal algebraic subgroup such that $\widehat G=\widehat D\cdot\widehat L$ where $\widehat L$ is the linear part of $\widehat G$.
Moreover, $\widehat D$ is connected and commutative.
Let $D=f^{-1}(\widehat D(\C))$.
So $G/D$ is definably isomorphic to $(\widehat G/\widehat D)(\C)$ and $G=D\cdot L$.
Arguing as we did in the proof of part~(a) one deduces the desired unique minimality properties of $D$ in $G$, and then that it is $B$-definable.
\end{proof}

\bigskip
\section{Abelian reduction}

\noindent
In what follows, by an {\em abelian variety} we mean simply a projective algebraic group.
That is, we do not assume connectedness.

\begin{definition}
Suppose $p\in S(B)$ is $\C$-internal.
We say that $p$ {\em has binding group an abelian variety in $\C$} to mean that $\aut(p/\C)$ is definably isomorphic
to $A(\C)$ for some abelian variety $A$ over $\C$.
\end{definition}

\begin{remark}
\label{connectedoverB}
A consequence of Corollary~\ref{connectedbind} is that if, in addition, $B=\acl(B)$ and $p$ is weakly $\C$-orthogonal then $A$ can be taken to be a connected abelian variety.
Moreover, because of the commutativity of the binding group, Lemma~\ref{overB} tells us that in this case we can actually take $A$ to be defined over $\dcl(B)\cap \C$ and for there to be a $B$-definable isomorphism between $\aut(p/\C)$ and $A(\C)$.
\end{remark}

Recall that, for us,  $\C$-internality includes stationarity.

\begin{lemma}
\label{abbyab}
Suppose $f:p\to q$ is a fibration of a $\C$-internal type $p\in S(B)$.
If both $q$ and the fibres of $f$ have binding groups abelian varieties in $\C$, then so does~$p$.
\end{lemma}

\begin{proof}
By Lemma~\ref{kernelstructure}, $\aut(p/\C)$ is definably an extension of $\aut(q/\C)$ by $\aut(r/\C)$, where $r:=\tp(\eta/Bf(\eta))$ and $\eta$ is a fundamental system for $p$ such that $f(\eta)$ is fundamental for $q$.
We are told that $\aut(q/\C)$ is an abelian variety in~$\C$.
As the extension of an abelian variety by an abelian variety is again an abelian variety, it suffices to observe that $\aut(r/\C)$ is an abelian variety in $\C$.

Writing $\eta=(a_1,\dots,a_n)$ and setting $r_i=\tp(a_i/Bf(a_i))$, we see that $\aut(r/\C)$ definably embeds in the cartesian product of the $\aut(r_i/\C)$.
But $r_i$ is a fibre of $f:p\to q$, and so by assumption $\aut(r_i/\C)$ is an abelian variety in $\C$.
As an algebraic subgroup of a cartesian product of abelian varieties is again an abelian variety, we deduce that $\aut(r/\C)$ is an abelian variety in $\C$, as desired.
\end{proof}

\begin{definition}
By an {\em abelian reduction} of a stationary type $p\in S(B)$ we mean a fibration $f:p\to p_{\ab}$ such that
\begin{enumerate}
\item
$p_{\ab}$ is $\C$-internal and has binding group an abelian variety in $\C$; and
\item
Universal Property: if $g:p\to q$ is any definable function with $q$ being $\C$-internal and having binding group an abelian variety in $\C$, then there is a definable function $h:p_{\ab}\to q$ such that $hf=g$.
\end{enumerate}
Sometimes we call $p_{\ab}$ itself an abelian reduction of $p$.
\end{definition}

It may be worth describing the universal property in terms of realisations: If $a\models p$ and $b\in\dcl(Ba)$ realises $p_{\ab}$, then the universal property is saying that whenever $c\in\dcl(Ba)$ is such that $\tp(c/B)$ is $\C$-internal with binding group an abelian variety in $\C$, then $c\in\dcl(Bb)$.

Note that an abelian reduction, if it exists, is unique up to interdefinability.

\begin{definition}
\label{trivar}
We say that a stationary type $p=\tp(a/B)$ has {\em trivial} abelian reduction if $p_{\ab}$ (exists and) is an algebraic type. (Since $p_{\ab}$ is stationary this would imply that it is realised in $\dcl(B)$.)
Equivalently, whenever $c\in\dcl(Ba)$ is such that $\tp(c/B)$ is $\C$-internal with binding group an abelian variety in $\C$, then $c\in\dcl(B)$.
\end{definition}

\begin{proposition}
\label{abfibab}
The fibres of the abelian reduction of a $\C$-internal type have trivial abelian reduction.
\end{proposition}

\begin{proof}
Suppose $p=\tp(a/B)$ is $\C$-internal and $b=f(a)\models p_{\ab}$.
The proposition is claiming that the (stationary) type $r:=\tp(a/Bb)$ has trivial abelian reduction.
Suppose $c\in\dcl(Ba)$ is such that $\tp(c/Bb)$ has binding group an abelian variety in $\C$.
Consider the $\C$-internal type $\tp(bc/B)$ and its fibration $g:\tp(bc/B)\to\tp(b/B)=p_{\ab}$ given by co-ordinate projection.
Both $p_{\ab}$ and the fibres of $g$ have binding group abelian varieties in $\C$, so that, 
by Lemma~\ref{abbyab}, so does $\tp(bc/B)$.
But $bc\in\dcl(Ba)$, and so, by the universal property of the abelian reduction, $c\in\dcl(Bb)$.
\end{proof}

\begin{theorem}
\label{abredexists}
Every $\C$-internal type admits an abelian reduction.
\end{theorem}

\begin{proof}
Fix $p=\tp(a/B)$ a $\C$-internal type.
Let
$\Ab_B(\C)$ denote the set of all finite tuples $b$ such that $\tp(b/B)$ is $\C$-internal and has binding group an abelian variety in $\C$.
Note that $\Ab_B(\C)$ is a $B$-invariant set.
Let $C:=\Ab_B(\C)\cap\dcl(Ba)$.

We first observe that if $a_1,\dots,a_n\in C$ then $(a_1,\dots,a_n)\in C$.
Indeed, note first of all that as $p$ is stationary and $a_1,\dots,a_n\in\dcl(Ba)$, we have that  $\tp(a_1,\dots,a_n/B)$ is stationary.
Moreover, it is $\C$-internal.
The binding group of $\tp(a_1,\dots,a_n/B)$ embeds definably into the cartesian product of the binding groups of the $\tp(a_i/B)$, each of which is an abelian variety over $\C$.
Hence, $(a_1,\dots,a_n)\in\Ab_B(\C)$ as desired.

Next we observe that $C$ is definably closed.
By the previous paragraph, it suffices to show that if $c\in C$ and $b=f(c)\in\dcl(c)$, then $b\in C$.
That $\tp(b/B)$ is (stationary and) $\C$-internal follows from the fact that $\tp(c/B)$ is.
Setting $r=\tp(c/B)$ and $q=\tp(b/B)$, we have a $B$-definable surjective homomorphism $f^\natural:\aut(r/\C)\to\aut(q/\C)$.
As the quotient of an abelian variety is an abelian variety, we conclude that $b\in\Ab_B(\C)$, as desired.

Next we observe that $C$ is relatively algebraically closed in $\dcl(Ba)$.
Suppose $b\in\acl(C)\cap\dcl(Ba)$.
Since $C$ is closed under finite tuples, there is $c\in C$ such that $b\in\acl(c)$.
Note that as $b,c\in\dcl(Ba)$, we have that $r=\tp(bc/B)$ is (stationary and) $\C$-internal.
Consider the finite-to-one definable function $f:r\to q=\tp(c/B)$ and the induced definable isogeny $f^\natural:\aut(r/\C)\to\aut(q/\C)$.
But $\aut(q/\C)$ is an abelian variety in $\C$, and any algebraic group isogeneous to an abelian variety is again an abelian variety.
Hence $\aut(r/\C)$ is an abelian variety in $\C$.
It follows that $(b,c)\in C$, and so $b\in C$ by definable-closedness.

We are now ready to prove the existence of an abelian reduction for $p=\tp(a/B)$.
Let $c$ be the canonical base of $\stp(a/C)$.
We claim that $\tp(c/B)=p_{\ab}$.
First of all, properties of the canonical base imply that $c\in\dcl(Ca)\cap\acl(C)$, see for example~\cite[1.2.26]{GST}.
That $C$ is relatively algebraically closed in $\dcl(Ba)$ thus implies that $c\in C$.
Hence $\tp(c/B)$ is $\C$-internal with binding group an abelian variety in $\C$.
Properties of the canonical base also imply that $\tp(a/Bc)$ is stationary.
So we have a fibration $f:p\to\tp(c/B)$.
It remains to show the universal property: given $b\in\dcl(Ba)$ with $\tp(b/B)$ having binding group an abelian variety over $\C$, we must show that $b\in\dcl(Bc)$.
But such $b$ is in $C$ by definition, and so, properties of the canonical base imply $a\ind_{Bc}b$.
That $b\in\dcl(Bc)$ then follows from $b\in\dcl(Ba)$ and stationarity.
\end{proof}

We end this section by considering how abelian reductions behave with respect to algebraic base extension.

\begin{lemma}
\label{passtoacl}
For any stationary type $q$ over $B$, let us denote by $\overline q$ the unique extension to $\acl(B)$.
If $p\in S(B)$ is $\C$-internal then $\overline p_{\ab}$ and $\overline{p_{\ab}}$ are interdefinable.
\end{lemma}

\begin{proof}
Here $\overline p_{\ab}$ is the abelian reduction of $\overline p$.
Let $p=\tp(a/B)$.
So $\overline p=\tp(a/\acl(B))$.
Since $p$ and $\overline p$ are $\C$-internal, $p_{\ab}$ and $\overline p_{\ab}$ exist by Theorem~\ref{abredexists}.
Let $p_{\ab}=\tp(b/B)$ and $\overline p_{\ab}=\tp(\overline b/\acl(B))$.
So $\overline{p_{\ab}}=\tp(b/\acl(B))$.
What the proposition is saying is that $b$ and $\overline b$ are interdefinable over $\acl(B)$.

One direction is relatively clear.
Note that $\overline{p_{\ab}}$ has binding group an abelian variety in $\C$ as, by Lemma~\ref{autbar}, $\aut(\overline{p_{\ab}}/\C)$ is of finite index in $\aut(p_{\ab}/\C)$.
The universal property for $\overline p_{\ab}$ thus implies that $b\in\dcl(\acl(B)\overline b)$.

To show that $\overline b\in\dcl(\acl(B) b)$, note first of all that $\overline b\in\dcl(\acl(B)a)\subseteq\acl(Ba)$ and hence has finitely many $Ba$-conjugates, say $c_1=\overline b,c_2,\dots c_n$.
Let $e$ be a code for the finite set $\{c_1,\dots,c_n\}$.
Then $e\in\dcl(Ba)$ and $\tp(e/B)$ is $\C$-internal.
Each $\tp(c_i/\acl(B))$, being conjugate to $\tp(\overline b/\acl(B))=\overline p_{\ab}$,  has binding group an abelian variety in~$\C$.
That is, in the notation of the proof of Theorem~\ref{abredexists}, $c_1,\dots,c_n\in\Ab_{\acl(B)}(\C)\cap\dcl(\acl(B)a)$.
That proof showed this set to be definably closed, so $e$, being in $\dcl(c_1\cdots c_n)$, is in $\Ab_{\acl(B)}(\C)$ too.
Hence $\tp(e/\acl(B))$ has binding group an abelian variety in $\C$.
But the binding group of $\tp(e/\acl(B))$ is of finite index in the binding group of $\tp(e/B)$ by Lemma~\ref{autbar}.
So $\tp(e/B)$ has binding group an abelian variety in $\C$.
By the universal property for $p_{\ab}$, it follows that $e\in\dcl(Bb)$.
Since $\overline b\in\acl(e)$, we get that $\overline b\in\acl(Bb)$.

So we have shown that $\overline b\in \acl(Bb)\cap \dcl(\acl(B)a)$.
It remains only to observe that $\acl(Bb)\cap \dcl(\acl(B)a)=\dcl(\acl(B) b)$.
Indeed, let $c \in \acl(Bb)\cap \dcl(\acl(B)a)$.
Since the abelian reduction is a fibration, $\tp(a/Bb)$ is stationary, and so $\tp(a/\acl(B)b)$ is stationary.
It follows that $\tp(c/\acl(B)b)$ is stationary.
But this type is algebraic, so that stationarity forces $c \in \dcl(\acl(B)b)$.
\end{proof}

\bigskip
\section{Binding group of the abelian reduction}

\noindent
Let us fix $B=\acl(B)$ and $p\in S(B)$ a $\C$-internal and weakly $\C$-orthogonal type.
We know by Corollary~\ref{connectedbind} that the binding group $G=\aut(p/\C)$ is connected, and we have the $B$-definable Chevalley decomposition
$$1\to L\to G\overset{\pi}\to A\to 0$$
given by Fact~\ref{definablechevalley}(a).
Now, let $f:p\to p_{\ab}$ be an abelian reduction given by Theorem~\ref{abredexists}.
So, like $A$, $\aut(p_{\ab}/\C)$ is definably isomorphic to the $\C$-points of a connected abelian variety over $\C$.
Our goal in this section is to prove that $\aut(p_{\ab}/\C)=A$.
That is, the binding group of the abelian reduction of $p$ is the abelian part of the binding group of $p$.

Here is a first approximation.

\begin{proposition}
\label{trivabred}
If $p$ has trivial abelian reduction then $G=L$.
\end{proposition}

\begin{proof}
Let $\eta$ be a fundamental system of realisations of $p$ given by Fact~\ref{makefundamental} so that $q:=\tp(\eta/B)$ is fundamental $\C$-internal, weakly $\C$-orthogonal, and $\aut(q/\C)=G$.
Writing $\eta=(a_1,\dots,a_n)$ we have by the Galois correspondence applied to $q$, namely Fact~\ref{gcorr}(a), that
$$H:=\{\sigma\in G:\text{ some (any) lifting of $\sigma$ to $\aut_B(\U/\C)$ fixes $a_1$}\}$$
is a $B\eta$-definable subgroup.
We consider $\pi(H)\leq A$.
Note that by rigidity of abelian varieties, and Remark~\ref{connectedoverB}, $\pi(H)$ is $B$-definable.

We first show that $\pi(H)=A$.
Since $L$ is normal in $G=\aut(q/\C)$, the Galois correspondence, namely Fact~\ref{gcorr}(b), gives us that $G/L=A=\aut(r/\C)$ where $r=\tp(b/B)$ is fundamental and $b\in \dcl(B\eta)$.
The Galois correspondence applied to $r$ implies that for some $c\in\dcl(Bb)$,
$$\pi(H)=\{\sigma\in A:\text{ some (any) lifting of $\sigma$ to $\aut_B(\U/\C)$ fixes $c$}\}$$
and $A/\pi(H)$ is the binding group of $\tp(c/B)$.
Note that 
$$H\subseteq\pi^{-1}\pi(H)=\{\sigma\in G:\text{ some (any) lifting of $\sigma$ to $\aut_B(\U/\C)$ fixes $c$}\}.$$
It follows that $c\in\dcl(Ba_1)$.
We therefore have a definable function $p\to \tp(c/B)$, and the latter has binding group an abelian variety in $\C$, namely $A/\pi(H)$.
The assumption that $p$ has trivial abelian reduction therefore implies that $c\in B$.
Which means that $\pi(H)=A$, as desired.

That $\pi(H)=A$ implies that $G=H\cdot L$.
It follows that the anti-linear part $D\leq G$ given by Fact~\ref{definablechevalley}(b)
must be contained in $H$.
Since $D$ is normal in $G$, we have $G/D=\aut(r/\C)$ where $r=\tp(b/B)$ is fundamental and $b\in\dcl(B\eta)$.
(We are repurposing the labels $r$ and $b$, these are not the same as they were in the previous paragraph.)
That $D\leq H$ implies $a_1 \in\dcl(Bb)$.
It follows that the binding group of $\tp(a_1/B)=p$, namely $G$, is a quotient of the binding group of $\tp(b/B)=r$, namely $G/D$.
This forces $D=1$ and $G$ to be definably isomorphic to the $\C$-points of a linear algebraic group over $\C$, as desired.
\end{proof}

We can now deduce:

\begin{theorem}
\label{abred=alb}
Suppose $B=\acl(B)$ and $p\in S(B)$ is $\C$-internal and weakly $\C$-orthogonal.
Then $\aut(p_{\ab}/\C)$ is the abelian part of $\aut(p/\C)$.
\end{theorem}

\begin{proof}
Let $f:p\to p_{\ab}$ be an abelian reduction.
We have the surjective definable homomorphism $f^\natural:\aut(p/\C)\to\aut(p_{\ab}/\C)$, with kernel $K$.
It suffices to show that $K$ is the linear part of $\aut(p/\C)$.
Since $\aut(p_{\ab}/\C)$ is definably isomorphic to the $\C$-points of an abelian variety over $\C$, we know that $K$ cannot be a proper subgroup of the linear part.
Hence, it suffices to show that $K$ is contained in the linear of $\aut(p/\C)$.

By Lemma~\ref{kernelstructure}, if we let $\eta\models p^{(n)}$ be fundamental for $p$ such that $f(\eta)$ is fundamental for $p_{\ab}$, and we let $r=\tp(\eta/Bf(\eta))$, then $K$ is definably isomorphic to $\aut(r/\C)$.
Writing $\eta=(a_1,\dots,a_n)$ and setting $r_i=\tp(a_i/Bf(a_i))$, we have a definable embedding of $\aut(r/\C)$ into the cartesian product of the $\aut(r_i/\C)$.
Now, each $r_i$ is a fibre of $f$, which by Proposition~\ref{abfibab} must therefore have trivial abelian reduction.
Let $\overline{r_i}$ denote the unique extension of $r_i$ to $\acl(Bf(a_i))$.
By Lemma~\ref{passtoacl}, $\overline{r_i}$ must also have trivial abelian reduction.
Note also that since $p$ is weakly $\C$-orthogonal so is $r_i$ and $\overline{r_i}$.
It follows by Proposition~\ref{trivabred} that each $\aut(\overline{r_i}/\C)$ is equal to its linear part.
As $\aut(\overline{r_i}/\C)$ is of finite index in $\aut(r_i/\C)$ by Lemma~\ref{autbar}, it follows that each $\aut(r_i/\C)$ is definably isomorphic to the $\C$-points of a linear algebraic group over $\C$.
The same is therefore true of $\aut(r/\C)$, and hence of $K$.
That is, $K$ is contained in the linear part of $\aut(p/\C)$.
\end{proof}

We end with a more or less  immediate application.
Recall from~\cite{moosa2014some} that a stationary type $p\in S(B)$ is said to {\em admit no proper fibrations} if whenever $p\to q$ is a fibration then either $p$ and $q$ are interdefinable or $q$ is realised in $\dcl(B)$.

\begin{remark}
Because we insist that the fibres of a fibration be stationary, our notion of admitting no proper fibrations is on the face of it weaker than~\cite[Definition~2.1]{moosa2014some}.
Nevertheless, the two definitions are equivalent: {\em If $p=\tp(a/B)$ admits no proper fibrations then whenever $c\in\dcl(Ba)\setminus\acl(B)$ we must have $a\in\acl(Bc)$.}
\end{remark}

\begin{proof}
Given such $c$, let $e$ be the canonical base of $\stp(a/Bc)$.
Then $e$ is contained in $\dcl(Ba)\cap\acl(Bc)$.
Letting $q:=\tp(e/B)$ we have a definable function $p\to q$ which is a fibration since $\tp(a/Be)$ is stationary.
Since $c\notin\acl(B)$, we must have that $a\nind_Bc$, and hence $e\notin\dcl(B)$.
So, as $p$ admits no proper fibrations in the sense we have defined above, we must have that $p$ and $q$ are interdefinable, from which it follows that $a\in\acl(Bc)$, as desired.
\end{proof}

In any case, Theorem~\ref{abred=alb} gives the following information about $\C$-internal types admitting no proper fibrations.

\begin{corollary}
\label{nofibre}
Suppose $p\in S(B)$ is $\C$-internal and admits no proper fibrations.
Then $\aut(p/\C)$ is equal to either its abelian or its linear part.
\end{corollary}

\begin{proof}
We may assume that $B=\acl(B)$.
Indeed, by Lemma~\ref{autbar}, if the binding group of the extension to $\acl(B)$ is either an abelian variety or a linear algebraic group in $\C$, then so is $\aut(p/\C)$.

Since $p$ admits no proper fibrations, letting $a\models p$, we have that either $a$ is a tuple of constants or $\dcl(Ba)\cap\C=B\cap\C$.
In the former case $\aut(p/\C)$ is trivial and we are done.
So we may assume that $\dcl(Ba)\cap\C=B\cap\C$.
In other words, $p$ is weakly $\C$-orthogonal.

Now, by Theorem~\ref{abredexists}, the fibration $p\to p_{\ab}$ exists, and by Theorem~\ref{abred=alb}, $\aut(p_{\ab}/\C)$ is the abelian part of $\aut(p/\C)$.
The assumption of no proper fibrations implies that either $p$ is interdefinable with $p_{\ab}$ or $p$ has trivial abelian reduction.
In the former case $\aut(p/\C)$ is equal to its abelian part, and in the latter case $\aut(p/\C)$ is equal to its linear part (by Proposition~\ref{trivabred}).
\end{proof}

\bigskip
\section{Base extension}

\noindent
We now consider how abelian reductions behave with respect to nonforking extensions.
We already dealt with extensions to the algebraic closure in Lemma~\ref{passtoacl}.
But to pass to arbitrary nonforking extensions we will need to use what we have learned from Theorem~\ref{abred=alb} about the binding group of abelian reductions.

For a stationary type $p\in S(B)$ and a set $E\supseteq B$, we denote by $p_E$ the unique nonforking extension of $p$ in $S(E)$.

\begin{proposition}
\label{abred-nf}
Suppose $B\subseteq E$ and $p\in S(B)$ is stationary.
If $p_E$ is $\C$-internal and weakly $\C$-orthogonal, then
$(p_{\ab})_E$ and $(p_E)_{\ab}$ are interdefinable.
\end{proposition}

\begin{proof}
By Lemma~\ref{passtoacl} we may assume that $B$ and $E$ are $\acl$-closed.
Indeed, assuming we have proved the result in that case, we have
\begin{eqnarray*}
(p_{\ab})_E
&=&\big((p_{\ab})_{\acl(B)}\big)_{\acl(E)}\upharpoonright_E\\
&=&\big((p_{\acl(B)})_{\ab}\big)_{\acl(E)}\upharpoonright_E\ \ \text{ by~Lemma~\ref{passtoacl}}\\
&=&\big((p_{\acl(B)})_{\acl(E)}\big)_{\ab}\upharpoonright_E\ \ \text{ by assumption}\\
&=&\big((p_E)_{\acl(E)}\big)_{\ab}\upharpoonright_E\\
&=&\big((p_E)_{\ab}\big)_{\acl(E)}\upharpoonright_E\ \ \text{ by~Lemma~\ref{passtoacl}}\\
&=&(p_E)_{\ab}.
\end{eqnarray*}
Here we are abusing notation as by equality in the above computations we sometimes (i.e., in lines $2,3$, and $5$) really only mean interdefinable.
In any case, we conclude that $(p_{\ab})_E$ and $(p_E)_{\ab}$ are interdefinable, as desired.

Note that since $p_E$ is $\C$-internal so is $p$.
In particular, the abelian reductions $f:p\to p_{\ab}$ and $g:p_E\to(p_E)_{\ab}$ exist.
The fibration $f$ lifts canonically to the nonforking extension so that we also have $f_E:p_E\to (p_{\ab})_E$.
Indeed, $f_E$ is just $f$ in the following sense:
If $a\models p$ with $a\ind_BE$ then $a\models p_E$, $f(a)\models (p_{\ab})_E$, and $\tp(a/Ef(a))$ is stationary as it is a nonforking extension of the stationary type $\tp(a/Bf(a))$.
On the other hand, we also have a definable function $h:(p_E)_{\ab}\to (p_{\ab})_E$ making the triangle commute:
$$\xymatrix{
& (p_E)_{\ab}\ar[dd]^h\\
p_E\ar[ur]^g\ar[dr]_{f_E}&\\
& (p_{\ab})_E
}$$
Indeed, this follows from the universal property for $g$ since the binding group of $(p_{\ab})_E$ embeds into the binding group of $p_{\ab}$ and hence is also an abelian variety in~$\C$.
We want to prove that $f_E:p_E\to (p_{\ab})_E$ is an abelian reduction of $p_E$.
This will force $h$ to be a definable bijection, witnessing the truth of the proposition.
But note that so far the only assumption we have used is that $p$ is $\C$-internal (and not yet that $p$, and even $p_E$, is weakly $\C$-orthogonal).

We pass to binding groups: the above triangle induces the following commuting triangle of surjective definable homomorphisms:
$$\xymatrix{
& A:=\aut((p_E)_{\ab}/\C)\ar[dd]^{h^\natural}\\
G:=\aut(p_E/\C)\ar[ur]^{g^\natural}\ar[dr]_{f_E^\natural}&\\
& A':=\aut((p_{\ab})_E/\C)
}$$
We claim now that $h^\natural:A\to A'$ is an isomorphism.
Since $E=\acl(E)$ and both $(p_E)_{\ab}$ and $(p_{\ab})_E$ are $\C$-internal and weakly $\C$-orthogonal (because $p_E$ is), these binding groups are connected by Lemma~\ref{connectedbind}.
It therefore suffices to show that $\dim A\leq \dim A'$.
Note that $G$ embeds into $\aut(p/\C)$.
Also, by Theorem~\ref{abred=alb}, $A$ is the abelian part of $G$.
Hence $A$ embeds in the abelian part of $\aut(p/\C)$, which is in turn equal to $\widetilde A:=\aut(p_{\ab}/\C)$.
So $\dim A\leq\dim\widetilde A$.
Now, $A'$ embeds in $\widetilde A$ as it is the binding group of a nonforking exctension of $p_{\ab}$.
As these groups are commutative, both homogeneous spaces $\big(A', (p_{\ab})_E(\U)\big)$ and $\big(\widetilde A, p_{\ab}(\U)\big)$ are principal.
It follows that $\dim A'=U((p_{\ab})_E)=U(p_{\ab})=\dim \widetilde A$.
This completes the proof that $h^\natural$ is an isomorphism.
In fact, we have exhibited canonical identifications $A=A'=\widetilde A$.
The upshot is that $A'$ is the abelian part of $G$ and so we have a definable Chevalley decomposition of the form
$$\xymatrix{
1\ar[r]&L\ar[r]&G\ar[r]^{f_E^\natural}&A'\ar[r]&0
}.$$

Let $a\models p_E$ and set $b:=f(a)\models (p_{\ab})_E$.
To show that $(p_{\ab})_E$ is an abelian reduction of $p_E$ we need to prove the universal property: given $c\in\dcl(Ea)$ such that $\tp(c/E)$ is $\C$-internal with binding group an abelian variety in $\C$, we need to show that $c\in\dcl(Eb)$.
Choose a fundamental system of realisations $\eta$ of $p_E$ extending $a$, and let $q:=\tp(\eta/E)$.
So $G=\aut(q/\C)$ also, and now we can apply the Galois correspondence.
The definably closed set $\dcl(Ec)$ corresponds to a $E\eta$-definable subgroup $H_c\leq G$.
But since $\tp(c/E)$ is fundamental (it has a commutative binding group), we know that in fact $H_c$ is $E$-definable and normal, and that $G/H_c$ is the binding group of $\tp(c/E)$.
In particular, $G/H_c$ is an abelian variety in $\C$, and so, from the definable Chevalley decomposition exhibited above, we must have that $L\leq H_c$.
On the other hand, the Galois correspondence will associate to $L=\ker(f_E^\natural)$ the definably closed set $\dcl(Eb)$, this is because it is $f_E$ that witnesses $b\in\dcl(E\eta)$.
Hence $L\leq H_c$ implies $c\in\dcl(Eb)$, as desired.
\end{proof}

\begin{remark}
It does not suffice in Proposition~\ref{abred-nf} to assume only that $p$ is weakly $\C$-orthogonal, we require the nonforking extension $p_E$ also be weakly $\C$-orthogonal.
Indeed, as long as $p$ has a proper abelian reduction, that is, as long as $p$ is not interdefinable with $p_{\ab}$, we must have that $p_E$ is not interdefinable with $(p_{\ab})_E$.
But if we take $E$ to be a model, for example, then $p_E$ is interdefinable with $(p_E)_{\ab}$.
\end{remark}

\bigskip
\section{Applications to differentially closed fields}
\label{sec:dcf}
\noindent
We specialise now to $T=\dcf$ and $\C$ the field of constants.
In fact, we could work more generally in the partial case of several commuting derivations (and $\C$ the total constants), but we stick to the ordinary case for notational convenience.

Our applications will focus on the {\em autonomous} case, that is,  when the parameters are contained in $\C$.
The main reason for this is the following important fact about binding groups in that case, which expresses the dramatic failure of the inverse Galois problem over the constants.
See, for example, the end of the proof of Proposition~4.9 in~\cite{JJP}, or the proof of Theorem~3.9 in~\cite{c3c2}.

\begin{fact}
\label{linearbg}
Suppose $k$ is an algebraically closed field of constants and $p\in S(k)$ is $\C$-internal and weakly $\C$-orthogonal.
If  $\aut(p/\C)$ is equal to its linear part then it is commutative.
\end{fact}

Note that if $p \in S(k)$ is $\C$-internal and weakly $\C$-orthogonal, the commutativity of $\aut(p/\C)$ implies that $p$ is fundamental as any faithful transitive commutative group action is principal. This observation will be used at several places in this section.
\medskip
\subsection{Bounding nonorthogonality}
Our work on abelian reductions has a surprising consequence for the structure of finite rank types over constant parameters:
every type nonorthogonal to the constants becomes non-weakly-orthogonal in the second Morley power.

\begin{theorem}
\label{nwodeg}
Suppose $B\subseteq\C$ and $p\in S(B)$ is stationary and of finite $U$-rank.
If $p$ is not $\C$-orthogonal then $p^{(2)}$ is not weakly $\C$-orthogonal.
\end{theorem}

\begin{proof}
We may assume that $p$ is weakly $\C$-orthogonal and that $B=\acl(B)$.
By nonorthogonality to $\C$ there is a definable function $p\to q$ with $q$ of positive $U$-rank and $\C$-internal, see~\cite[7.4.6]{GST}.
Now, there exists a fibration $q\to r$ such that $r$ is of positive $U$-rank and admits no proper fibrations.
Indeed, if $q$ itself admits no proper fibration then $r=q$ works, otherwise take a proper fibration $q\to q'$ and proceed by induction on $U$-rank, noting that $U(q')<U(q)$ (and that $U$-rank~$1$ types admit no proper fibrations).
So we have a definable function $p\to r$ with $r$ internal to~$\C$, weakly orthogonal to $\C$, and admitting no proper fibrations.
By Corollary~\ref{nofibre}, we conclude that $\aut(r/\C)$, which is connected by Corollary~\ref{connectedbind}, agrees with either its abelian part or its linear part.
By Fact~\ref{linearbg}, since we are working over constant parameters, if the binding group is linear then it is commutative.
Therefore, in either case, whether $\aut(r/\C)$ agrees with its abelian part or its linear part, it is a commutative group.
But commutativity implies that $(\aut(r/\C),r(\U))$ is a principal homogeneous space, and hence $r$ is fundamental.
In particular, $r^{(2)}$ is not weakly $\C$-orthogonal.
Hence $p^{(2)}$ is not weakly $\C$-orthogonal, as desired.
\end{proof}

The assumption that $p$ be autonomous is necessary.
For example, for each $n\geq 2$, a (nonautonomous) $\C$-internal type $p$ exists whose binding group action is definably isomorphic to the natural action of $\PSL_{n+1}(\C)$ on $\PP^n(\C)$.
See, for example, Section~4.2 of~\cite{nmdeg} for a construction.
Now, that action is {\em generically $(n+2)$-transitive}, meaning that the co-ordinate-wise action on the $(n+2)$nd cartesian power of $\PP^n(\C)$ admits a Zariski open orbit.
But this implies that $\aut(p/\C)$ acts transitively on $p^{(n+2)}(\U)$, and hence $p^{(n+2)}$ is weakly $\C$-orthogonal.
So, in the general nonautonomous situation, one cannot expect a bound independent of dimension on how high a Morley power must be considered to verify $\C$-orthogonality.

Let us formulate Theorem~\ref{nwodeg} geometrically.
Fix an algebraically closed subfield $k\subseteq\C$.
Recall that every finite dimensional type in $\dcf$ over $k$ arises, up to interdefinability, as the generic type of an algebraic vector field $(X,v)$ over~$k$.
Here $X$ is an irreducible affine algebraic variety over $k$ and $v:X\to TX$ is a regular section to the tangent bundle.
The {\em generic type} of $(X,v)$ over $k$ is the type $p(x)$ saying that $x\in X$, that $x$ is not contained in any proper subvariety of $X$ over~$k$, and that $v(x)=(x,\delta(x))$.
This type is of finite $U$-rank, bounded above by $\dim X$.
Weak $\C$-orthogonality of $p$ translates to the fact that $(X,v)$ admits no nontrivial {\em rational first integrals}, i.e., there is no nonconstant rational function $f$ on $X$ with $\delta_v(f)=0$.
Here $\delta_v$ is the $k$-linear derivation on $k(X)$ induced by the vector field~$v$.
Outright $\C$-orthogonality of $p$ corresponds to the condition that, for any other algebraic vector field $(Y,w)$, every rational first integral of the cartesian product $(X\times Y,v\times w)$ lifts a rational first integral of $(Y,w)$.

\begin{corollary}
\label{nwodeg-geometric}
Suppose $(X,v)$ is an irreducible algebraic vector field over an algebraically closed field $k$ of characteristic zero.
The following are equivalent:
\begin{itemize}
\item[(i)]
$(X^2,v^2)$ admits no nontrivial rational first integral,
\item[(ii)]
$(X^n,v^n)$ admits no nontrivial rational first integral, for all $n\geq 1$,
\item[(iii)]
for all irreducible algebraic vector fields $(Y,w)$ over $k$ admitting no nontrivial rational first integrals, the cartesian product $(X\times Y,v\times w)$ admits no nontrivial rational first integral.
\end{itemize}
\end{corollary}

\begin{proof}
Clearly~(ii) implies~(i).
To see that~(iii) implies~(ii), apply~(iii) first to $Y$ a one-point variety with the trivial vector field structure to see that $(X,v)$ itself must admit no nontrivial rational first integrals.
Now apply~(iii) to $(Y,w)=(X,v)$ and deduce~(ii) by induction on $n$.

So it remains to show that~(i) implies~(iii).
For this, we work in a sufficiently saturated model $\U\models\dcf$ with field of constants $\C$ extending~$k$ and consider $p\in S(k)$ the generic type of $(X,v)$.
Condition~(i) says that $p^{(2)}$ is weakly $\C$-orthogonal, which by Theorem~\ref{nwodeg} implies that $p$ is $\C$-orthogonal.
As explained above, this means that every rational first integral of $(X\times Y,v\times w)$ comes from a rational first integral of $(Y,w)$, and since we are assuming there are no nontrivial such, we deduce~(iii).
\end{proof}

\medskip
\subsection{Classifying two-dimensional types}

\noindent
In this section we deduce, from our results on abelian reductions, a classification of $2$-dimensional autonomous types that are $\C$-internal and weakly $\C$-orthogonal.
In the fundamental $\C$-internal case such a classification in terms of logarithmic-differential equations on algebraic groups is a known consequence of differential Galois theory, even in arbitrary dimension and over possibly nonconstant parameters.

Let us briefly recall logarithmic-differential equations on algebraic groups, see~\cite{pillay2004} for details.
If $G$ is a connected algebraic group over the constants then we have a {\em logarithmic-derivative on~$G$} which is a certain $\delta$-rational surjective crossed homomorphism
$\dlog_G:G(\U)\to LG(\U)$
to the Lie algebra of $G$.
The fibre of $\dlog_G$ above~$0$ is $G(\C)$ and all the other fibres are left cosets of $G(\C)$.
By a {\em logarithmic-differential equation on~$G$} we mean an equation of the form $\dlog_G(x)=\rho$
where $\rho\in LG(\U)$.
The generic type, $p$, of this equation is therefore $\C$-internal.
The equation is said to be {\em full} if in addition $p$ is weakly $\C$-orthogonal.
In that case $p$ is isolated by $\dlog_G(x)=\rho$, and the action of $\aut(p/\C)$ is definably isomorphic to the action of $G(\C)$ by multiplication on the right.
The following model-theoretic articulation of a theorem of Kolchin's~\cite[$\S$VI.9]{KolchinDAAG} says that every fundamental $\C$-internal weakly $\C$-orthogonal type arises this way.

\begin{fact}[Kolchin]
\label{alldlog}
Suppose $k$ is an algebraically closed differential field and $p$ is fundamental $\C$-internal and weakly $\C$-orthogonal complete type over $k$.
Then there is a connected algebraic group $G$ over $k\cap\C$ such that $p$ is interdefinable with the generic type of a full logarithmic-differential equation on $G$ over $k$.
\end{fact}

Our goal here is to improve on this classification result, when we are in dimension~$2$ and over constant parameters, by dropping the assumption that the type be fundamental.
The idea is to consider the abelian reduction $f: p\to p_{\ab}$ and thereby reduce to the fundamental $p_{\ab}$ (which can be handled by Fact~\ref{alldlog}) along with the fibre of $f$ which in the proper nontrivial case will be of dimension~$1$.
So we will require a complete understanding of dimension~$1$ (possibly nonfundamental) $\C$-internal and weakly $\C$-orthogonal types over possibly nonconstant parameters.
We therefore begin with such an analysis.

Here the central examples are {\em Riccati equations}, order~1 equations in a single variable of the form $x'=ax^2+bx+c$.
We can, and do, always write such an equation as $x' = -\gamma x^2 + (\alpha - \delta)x + \beta$,
where $M = \begin{pmatrix} \alpha & \beta \\ \gamma & \delta \end{pmatrix}\in\SL_2(k)$.
The reason for doing so is that this is the equation satisfied by $\frac{y_1}{y_2}$ where $Y=\begin{pmatrix} y_1 \\ y_2\end{pmatrix}$ is a nonzero solution to the system of linear order~1 differential equations given by $Y'=MY$.
In any case, almost all $1$-dimensional $\C$-internal weakly $\C$-orthogonal types come from Riccati equations:

\begin{proposition}
\label{classification-order one}
Let $k$ be an algebraically closed differential field.
Suppose $p$ is a complete $1$-dimensional type over $k$ that is $\mathcal C$-internal, weakly $\mathcal C$-orthogonal, and with $\aut(p/\C)$ equal to its linear part.
Then $p$ is interdefinable with the generic type of a Riccati equation.
\end{proposition}

\begin{proof}
Note that $G:=\aut(p/\C)$ is connected by Corollary~\ref{connectedbind}.
Since $p$ is an isolated type of $U$-rank $1$, $S:=p(\U)$ is a strongly minimal set that is a definable homogeneous space for $G$.
By Hrushovski's~\cite{hrushovski1989almost} classification of strongly minimal homogeneous spaces, we have that $(G,S)$ is definably isomorphic to one of the five possibilities:
\begin{enumerate}
\item\label{elliptic}
the action of an elliptic curve on itself by translation,
\item\label{Ga}
the action of $\Ga(\C)$ on $\C$ by translation,
\item\label{Gm}
the action of $\Gm(\C)$ on $\C\setminus\{0\}$ by multiplication,
\item\label{Aff}
the action of $\Aff_2(\C)$ on $\C$ by affine transformations, or
\item\label{PSL}
the action of $\PSL_2(\C)$ on $\PP(\C)$ by linear fractional transformations.
\end{enumerate}
The first case is ruled out by our assumption that $G$ is linear.
In cases~(\ref{Ga}) and~(\ref{Gm}) the binding group is commutative, and so $p$ is in addition fundamental, and hence Fact~\ref{alldlog} applies.
This gives that $p$ is interdefinable with  the generic type of the Riccati equation $x' = \beta$ in case~(\ref{Ga}) and $x' = \lambda x$ in case~(\ref{Gm}).

Let us consider case~(\ref{Aff}) and denote by $\phi$ a formula isolating $p$.
Since $\Aff_2(\C)$ acts uniquely $2$-transitively on $\C$, the type $p^{(2)}$ is isolated by
$$\psi(x,y) := \phi(x) \cup \phi(y) \cup \lbrace x \neq y \rbrace$$
and the action of $G$ on $p^{(2)}$ is principal. Indeed, by unique $2$-transitivity of the action of $G$ on $\phi(\U)$, the action of $G$ on $\psi(\U)$ is principal and since this is an action by automorphisms over $B\mathcal C$, all the elements of $\psi(\U)$ realize the same type over $B\mathcal C$ and hence the type $p^{(2)}$ over $B$. It follows that $p^{(2)}$ is $\C$-internal, weakly $\C$-orthogonal and fundamental, with the same binding group as $p$
and we can apply Fact~\ref{alldlog} to it.
So $q$ is interdefinable with the generic type $r\in S(k)$ of a full logarithmic-differential equation on $\Aff_2$.
We view $\Aff_2$ as embedded in $\SL_2$ by representing an affine linear transformation $z\mapsto \alpha z+\beta$, where $\alpha\neq 0$, as the special linear matrix  $\begin{pmatrix} \alpha & \beta \\ 0 & \frac{1}{\alpha} \end{pmatrix}$.
Fix $b= \begin{pmatrix} b_1 & b_2 \\ 0 & \frac{1}{b_1} \end{pmatrix}\models r$.
Then $\begin{pmatrix} b_2 \\  \frac{1}{b_1} \end{pmatrix}$ satisfies an equation of the form $Y'=MY$ and hence $\xi:= b_1b_2$ satisfies a Riccati equation.
We will show that $p$ is interdefinable with $\tp(\xi/k)$.

As $r$ is interdefinable with $q=p^{(2)}$, we have $a\models p$ such that $a\in\dcl(kb)$.
By the Galois correspondence, $\dcl(ka)$ corresponds to a $1$-dimensional $kb$-definable subgroup $H_a\leq\aut(r/\C)$.
On the other hand, $\dcl(k\xi)$ also gives rise to a $1$-dimensional $kb$-definable subgroup $H_\xi\leq\aut(r/\C)$.
To show that $a$ and $\xi$ are interdefinable over $k$ it suffices to show that $H_a=H_\xi$.

First of all, $H_a$ is definably isomorphic to $\Gm(\C)$ in $\Aff_2(\C)$.
Indeed, $H_a$ is, under the natural identification of $\aut(r/\C)$ with $G$, the stabiliser of $a$, and in case~(\ref{Aff}) the stabilisers of all points are of this form.
Since all the copies of $\Gm$ in $\Aff_2$ are conjugate, and conjugation amounts to replacing $b$ by another realisation of~$r$, we may assume that if $\sigma\in H_a$ then $\sigma$ acts on the matrix $b$ by multiplication on the right by a matrix of the form $\begin{pmatrix} c & 0 \\ 0 & \frac{1}{c}\end{pmatrix}$ for some nonzero $c\in\C$.
Choosing a lifting $\widehat\sigma\in\aut_k(\U/\C)$, we can now compute $\widehat\sigma(\xi)$ as follows:
From
$$  \begin{pmatrix} b_1 & b_2 \\ 0 & \frac{1}{b_1} \end{pmatrix} \begin{pmatrix} c & 0 \\ 0 & \frac{1}{c} \end{pmatrix} = \begin{pmatrix} cb_1 & \frac{b_2}{c} \\ 0 & \frac{1}{cb_1} \end{pmatrix}$$
we see that $\widehat\sigma(b_1)=cb_1$ and $\widehat\sigma(b_2)=\frac{b_2}{c}$.
So $\widehat\sigma(\xi)=\widehat\sigma(b_1b_2)=b_1b_2=\xi$.
What we have shown is that  $H_a\leq H_\xi$.
But in $\Aff_2$ the only $1$-dimensional algebraic subgroup containing $\Gm$ is $\Gm$ itself.
So $H_a=H_\xi$, as desired.

The final case~(\ref{PSL}) is handled similarly. This time, since $\PSL_2(\C)$ acts uniquely $3$-transitively on $\PP(\C)$, a similar reasoning as in case (4) shows that $q:=p^{(3)}$ is weakly $ \mathcal C$-orthogonal and fundamental and so interdefinable with $r\in S(k)$ the generic type of a full logarithmic-differential equation on $\PSL_2$.
We represent the elements of $\PSL_2$ as equivalence classes $\begin{bmatrix} \alpha & \beta \\ \gamma & \delta \end{bmatrix}$
of $\begin{pmatrix} \alpha & \beta \\ \gamma & \delta \end{pmatrix}\in\operatorname{GL}_2$ under the action of scalar multiplication.
Fixing a realisation $b=\begin{bmatrix} b_1 & b_2 \\ b_3 & b_4 \end{bmatrix}$ of~$r$, we consider $\xi:=\frac{b_1}{b_3}$, which again satisfies a Riccati equation, and show that $p$ is interdefinable with $\tp(\xi/k)$.

We have $a\models p$ with $a\in\dcl(kb)$ giving rise via Galois correspondence to a $2$-dimensional $kb$-definable subgroup $H_a\leq\aut(r/\C)$, and $\xi\in\dcl(kb)$ giving rise to a $2$-dimensional $kb$-definable subgroup $H_\xi\leq\aut(r/\C)$.
We aim to show that $H_a=H_\xi$.
As all $2$-dimensional algebraic subgroups of $\PSL_2$ are isomorphic to $\Aff_2$, and all such are conjugate, after possibly changing the realisation $b\models r$ we may assume that if $\sigma\in H_a$ then $\sigma$ acts on $b$ by multiplication on the right by $ \begin{bmatrix} c_1 & c_2 \\ 0 & \frac{1}{c_1} \end{bmatrix}\in\Aff_2(\C)$.
Since
$  \begin{bmatrix} b_1 & b_2 \\ b_3 & b_4 \end{bmatrix} \begin{bmatrix} c_1 & c_2 \\ 0 & \frac{1}{c_1} \end{bmatrix} = \begin{bmatrix} c_1b_1 &c_2b_1 + \frac{b_2}{c_1}  \\   c_1b_3 & c_2b_3 + \frac{b_4}{c_1}\end{bmatrix}$
we see that $\sigma$, or rather any lifting of it to $\aut_k(\U/\C)$, must fix $\frac{b_1}{b_3}=\xi$.
Hence, $H_a\leq H_\xi$.
As $H_\xi$ is connected and also $2$-dimensional, being isomorphic to $\Aff_2(\C)$, we must have $H_a=H_\xi$, as desired.
\end{proof}

Combining this understanding of the $1$-dimensional case with our results on abelian reductions, we obtain the following classification of the autonomous $2$-dimensional case.

\begin{theorem}
\label{dim2}
Suppose $k$ is an algebraically closed subfield of the constants and $p\in S(k)$ is a $2$-dimensional $\C$-internal weakly $\C$-orthogonal type.
Then one of the following two cases occurs:
\begin{enumerate}
\item
$p$ is interdefinable with the generic type of a logarithmic-differential equation on either $\mathbb G_m^2$, $\mathbb G_m\times\mathbb G_a$, or a connected abelian surface; or
\item
there is a fibration $f:p\to q$ where $q$ is the generic type of a logarithmic-differential equation on an elliptic curve and, for $a\models p$, $\stp(a/kf(a))$ is interdefinable with the generic type of a Riccati equation.
\end{enumerate}
\end{theorem}

\begin{proof}
First note that $\aut(p/\C)$ is connected by Corollary~\ref{connectedbind}.
By Theorems~\ref{abredexists} and~\ref{abred=alb} we have that an abelian reduction $f:p\to p_{\ab}$ exists and $\aut(p_{\ab}/\C)$ is the abelian part of $\aut(p/\C)$.
The dimension of $p_{\ab}$ is $0,1$, or $2$, and we analyse each case separately.

If $p_{\ab}$ is trivial then $\aut(p/\C)$ is equal to its linear part, and hence commutative by Fact~\ref{linearbg}.
It follows that $p$ is fundamental and Fact~\ref{alldlog} applies.
So $p$ is interdefinable with the generic type of a logarithmic-differential equation on $\mathbb G_a^2, \mathbb G_m^2$ or $\mathbb G_m\times\mathbb G_a$.
The first of these cannot occur as the generic type of a logarithmic-differential equation on $\mathbb G_a^2$ is never weakly $\C$-orthogonal.

If $p_{\ab}$ is of dimension $2$ then we may assume that $p=p_{\ab}$.
In that case the binding group is definably isomorphic to the $\C$-points of a connected abelian variety $A$ over $k\cap \C$.
Again, $p$ is fundamental by commutativity of the binding group, and Fact~\ref{alldlog} tells us that $p$ is definably isomorphic to the generic type of a logarithmic-differential equation on $A$.

Finally, suppose $\dim(p_{\ab})=1$.
Then, as $p_{\ab}$ is fundamental (having commutative binding group), it is interdefinable with the generic type $q$ of a logarithmic-differential equation on an elliptic curve.
We may take $p_{\ab}=q$ and consider the fibration $f:p\to p_{\ab}$.
Letting $a\models p$, the fibre $\tp(a/kf(a))$ is a $1$-dimensional stationary type that is $\C$-internal, weakly $\C$-orthogonal, and of trivial abelian reduction.
Letting $K=\acl(kf(a))$, Proposition~\ref{classification-order one} therefore applies to $\tp(a/K)=\stp(a/kf(a))$, which is therefore interdefinable with the generic type of a Riccati equation.
\end{proof}

Let us express this theorem in terms of the birational geometry of vector fields.
Recall that a finite dimensional type $p$ over constant parameters is (up to interdefinability) the generic type of an algebraic vector field $(X,v)$.
As we have already mentioned, that $p$ is weakly $\C$-orthogonal means that $(X,v)$ admits no nontrivial rational first integrals.
That $p$ is $\C$-internal says that $(X,v)$ is {\em isotrivial}: that is after some base extension to a possibly non constant differential field extension $K$ of $k$, $(X,v)$ is birationally equivalent to the base change to $K$ of a trivial vector field defined over the constant field of $K$,  namely where the section is the zero-section of the tangent bundle.   
The condition of isotriviality appears in  seemingly different but equivalent forms in the literature.
In an analytic setting,  Umemura defines the class of differential equations whose {\em general solution depends rationally on the initial conditions} (see~\cite[Definition~7]{umemura}) which coincides in the autonomous case with the class of isotrivial algebraic vector field; see Condition~(9) on p.~134 of~\cite{umemura}.
In a differential Galois-theoretic setting,  the isotrivial algebraic vector fields are also characterised among the algebraic vector fields without nontrivial first rational integrals as those which admit a generic solution in a strongly normal extension in the sense of Kolchin.

Theorem~\ref{dim2} is precisely about $2$-dimensional algebraic vector fields that are isotrivial and admit no nontrivial rational first integrals; it gives a birational classification of these:

\begin{corollary}
\label{classification-ordertwo}
Suppose $k$ is an algebraically closed field of characteristic zero and $(X,v)$ is an irreducible algebraic vector field over $k$ of dimension $2$ that is isotrivial and admits no nontrivial rational first integrals.
Then one of the following holds:
\begin{enumerate}
 \item
 $(X,v)$ is birationally equivalent over $k$ to
 \begin{itemize}
 \item[(a)]
 the vector field on $\mathbb A^2$ given by $w(x,y)=(x,y,\lambda x, \mu y)$ for some $\mathbb Z$-linearly independent  $\lambda,\mu\in k$, or
 \item[(b)]
 the vector field on $\mathbb A^2$ given by $w(x,y)=(x,y,\lambda, \mu y)$ for some nonzero $\lambda,\mu\in k$, or,
\item[(c)]
a translation invariant vector field on a connected abelian surface, or
\end{itemize}
\item
$X$ admits a dominant rational map $X\to E$ to an elliptic curve $E$ over $k$ such that $v$ lifts a translation invariant vector field $w$ on $E$, and such that after base extension to an \'etale open set $U\to E$, $(X,v)$ becomes birationally equivalent to 
the projectivisation of a vector field on $\AA^2 \times U$ that lifts $w_U$ and is linear on the fibres over $U$.
\end{enumerate}
\end{corollary}

Before giving the proof, let us give some explanation on the terminology used.
First, by a {\em translation invariant} vector field on a commutative algebraic group $G$ is meant a vector field $v:G\to TG$ such that $v\circ \lambda=(d\lambda)\circ v$ for any translation $\lambda:G\to G$.
Such vector fields are of course determined by $\rho=v(0)\in LG$, and they correspond precisely to the logarithmic-differential equation $\dlog_G(x)=\rho$ in the sense that the set of $D$-points of $(G,v)$ is precisely the set of solutions to $\dlog_G(x)=\rho$.
Regarding projectivisation, if $w$ is an algebraic vector field on $\AA^{n+1}\times X$, for some variety $X$, that lifts a vector field $w_X$ on $X$ and is linear on the fibres, then projectivisation induces a vector field on $\PP^n\times X$ also lifting~$w_X$, which we denote by $(\PP^n\times X,\PP(w))$.

\begin{proof}[Proof of Corollary~\ref{classification-ordertwo}]
We work in a saturated model $\U\models \dcf$ with field of constants $\C$ extending $k$, and we
let $p\in S(k)$ be the generic type of $(X,v)$ over $k$.
The assumptions on $(X,v)$ imply that $p$ is $\C$-internal and weakly $\C$-orthogonal.
So Theorem~\ref{dim2} applies.

If we are in case~(1) of Theorem~\ref{dim2} then $p$ is interdefinable with the generic type of a logarithmic-differential equation on an algebraic group $G$ that is either $\mathbb G_m^2, \mathbb G_m\times\mathbb G_a$, or a connected abelian surface.
Hence $(X,v)$ is birationally equivalent to $(G,w)$ with $w$ an invariant vector field on $G$.
This gives rise to cases~(1)(a), (1)(b), and~(1)(c) of the Corollary, respectively.

Suppose, therefore that we are in case~(2) of Theorem~\ref{dim2}.
We have a fibration $f:p\to q$ where $q$ is the generic type of a logarithmic-differential equation on an elliptic curve $E$ over $k$, so of $(E,w)$ where $w$ is a translation invariant vector field on $E$ over $k$.
The fibration $f$ induces a dominant rational map $X\to E$ such that $v$ lifts $w$.
Moreover, for $a\models p$ and letting $b=f(a)$, we have that $\stp(a/kb)$ is interdefinable with the generic type of a Riccati equation over $\acl(kb)$.
So there is a finite extension $K$ of $k(b)$ such that $\tp(a/K)$ is interdefinable with the generic type $r$ of the Riccati equation corresponding to some $M_b\in\SL_2(K)$.
The extension $K$ corresponds to an \'etale open set $U\to E$, and  $w$ lifts uniquely to a vector field $w_U$ on $U$.
Shrinking~$U$ further if necessary, $w_U$ lifts to a vector field $\theta$ on $\mathbb A^2\times U$ that is linear on the fibres over $U$ and agrees with $M_b$ on the generic fibre.
The interdefinability of $\tp(a/K)$ and $r$ induces a birational equivalence between the base extension of $(X,v)$ to $(U,w_U)\to (E,w)$ and $(\PP^1\times U,\PP(\theta))$.
\end{proof}

\medskip
\subsection{No fibrations}
Our third application is more or less a direct combination of Corollary~\ref{nofibre} with Facts~\ref{linearbg} and~\ref{alldlog}.
It has to do with types $p\in S(B)$ that not only admit no proper fibrations themselves, but whenever $q\to p$ is any finite-to-one definable function then $q$ also admits no proper fibrations.
Put another way: for $a\models p$, if $c\in\acl(Ba)\setminus\acl(B)$ then $a\in \acl(Bc)$.
We say that $p$ {\em admits no proper almost-fibrations}.

\begin{theorem}
\label{nofibre-dcf}
Suppose $k\subseteq\C$ is an algebraically closed field and $p\in S(k)$ is stationary and nonalgebraic.
If $p$ admits no proper almost-fibrations then $p$ is either minimal or admits a finite-to-one definable function $f:p\to q$ where $q$ is the generic type of a logarithmic-differential equation on a simple abelian variety of dimension~$>1$.
\end{theorem}

\begin{proof}
It was shown in~\cite[Proposition~2.3]{moosa2014some} in an abstract finite rank setting that if $p$ admits no proper fibrations then either
\begin{itemize}
\item[(i)]
$p$ is interalgebraic with $r^{(k)}$ for some locally modular minimal type $r\in S(k)$, or
\item[(ii)]
$p$ is almost internal to a non locally modular minimal type.
\end{itemize}
Because we are assuming that $p$ admits no almost-fibrations, in case~(i) we would have that $r^{(k)}$ admits no proper fibrations as well, and hence $k=1$.
But then $p$ is minimal, as desired.
So we are in case~(ii).
This means that $p$ is almost $\C$-internal.
It follows that there is a $\C$-internal $q\in S(k)$ and a finite-to-one definable function $p\to q$.
Note that $q$ admits no proper fibrations either, and hence in particular is weakly $\C$-orthogonal.
Moreover, by Corollary~\ref{nofibre}, $\aut(q/\C)$ is equal to either its abelian or linear part.
By Fact~\ref{linearbg}, we have in either case that $\aut(q/\C)$ is commutative, and so $q$ is fundamental.
So, Fact~\ref{alldlog} applies, and $q$ is interdefinable with the generic type of a logarithmic-differential equation on a connected algebraic group $G$ over $k$.
Note that $G(\C)=\aut(q/\C)$, and hence $G$ is commutative.
On the other hand, no proper fibrations forces $G$ to be a simple algebraic group; the quotient modulo any proper infinite algebraic subgroup would give rise to a proper fibration.
So $G$ is either the additive group, the multiplicative group, or a simple abelian variety.
In the first two cases $q$, and hence $p$, is minimal.
This leaves the third possibility, namely that $q$ is interdefinable with the generic type of a logarithmic-differential equation on a simple abelian variety.
\end{proof}

\begin{remark}
It is worth noting that when $\dim(p)>1$ then in the minimal case of the dichotomy given by Theorem~\ref{nofibre-dcf} the $\acl$-geometry on the minimal type $p$ must be trivial.
This follows from the work of Hrushovski and Sokolovic~\cite{HrSo} on minimal types in $\dcf$.
\end{remark}

Again the assumption that $p$ be autonomous is necessary.
The examples constructed in~\cite[Section~4.2]{nmdeg} of (nonautonomous) $\C$-internal types $p$ whose binding group action are definably isomorphic to the natural action of $\PSL_{n+1}(\C)$ on $\PP^n(\C)$ are neither minimal (they are of rank~$n$) nor finite covers of logarithmic-differential equations on abelian varieties.
And they do not admit any proper almost-fibrations by the work in~\cite{nmdeg}, see Lemma~3.2 and Section~4.2 of that paper.

Here is a geometric formulation of Theorem~\ref{nofibre-dcf} whose verification we leave to the reader.
By a {\em generically finite cover} of an irreducible algebraic vector field $(X,v)$, we mean another irreducible algebraic vector field, $(Y,w)$, together with a dominant morphism $(Y,w)\to(X,v)$ whose generic fibres are finite (equivalently $\dim Y=\dim X$).

\begin{corollary}
\label{nofibre-dcf-geometric}
Suppose $(X,v)$ is an irreducible algebraic vector field over an algebraically closed field $k$ of characteristic zero.
Assume that for every generically finite cover $(Y,w)\to (X,v)$ over $k$, $(Y,w)$ admits no dominant morphism $(Y,w)\to (Z,u)$ over~$k$ with $0<\dim Z<\dim Y$.
Then either:
\begin{itemize}
\item[(i)]
for all irreducible algebraic vector fields $(Y,w)$ over $k$, every proper invariant subvariety of $(X\times Y, v\times w)$ over $k$ that projects dominantly onto both $X$ and $Y$ is a generically finite cover of $(Y,w)$, or
\item[(ii)]
$(X,v)$ is a generically finite cover of a translation invariant vector field on a simple abelian variety over $k$ of dimension $>1$.
\end{itemize}
\end{corollary}

\medskip
\subsection{Projective varieties admitting vector fields without first integrals}
Let us fix an irreducible projective algebraic variety $X$ over an algebraically closed field $k$ of characteristic zero.
We know that, because of projectivity, every algebraic vector field on $X$ is isotrivial; this is a theorem of Buium~\cite{Buium}.
Nevertheless, it is still possible for $(X,v)$ to admit no nontrivial rational first integrals -- that is, for its generic type in $\dcf$ to be weakly orthogonal to the constants.
For example, every nontrivial translation invariant vector field on a simple abelian variety is of this form.
It is natural to ask, therefore, what constraints the existence of such a vector field puts on $X$.
We give an answer here: $X$ must look very much like a projective vector bundle over its Albanese.

Recall that the {\em Albanese morphism of $X$}, with respect to a fixed $x\in X(k)$,  is a morphism $\alpha:X\to \Alb(X)$ where $\Alb(X)$ is a connected abelian variety over~$k$, $\alpha(x)=0$, and any other such morphism from $X$ to a connected abelian variety factors uniquely through $\alpha$.

The following is inspired by an analogous result for compact K\"ahler manifolds due to Fujiki~\cite[Theorem~1]{fujiki83}.

\begin{theorem}
\label{albanese-surjective}
Suppose $X$ is a smooth  projective variety over an algebraically closed field $k$ of characteristic zero.
If $X$ admits an algebraic vector field with no nontrivial rational first integrals then the Albanese morphism $\alpha:X\to\Alb(X)$ is surjective with generic fibre absolutely irreducible and birationally equivalent to an algebraic homogeneous space for a connected linear algebraic group.\footnote{Following Fujiki~\cite{fujiki83} we might say the generic fibre is {\em almost homogeneous unirational}.}
\end{theorem}

\begin{proof}
First of all, we claim that it suffices to show that there is some surjective morphism $f:X\to A$ to a connected abelian variety, with generic fibre absolutely irreducible and birationally equivalent to an algebraic homogeneous space for a connected linear algebraic group.
Indeed, such an $f$ would have to be (up to isomorphism) the Albanese morphism.
Let us give some details:
Since $f$ is surjective, fixing $x\in A$ such that $f(x)=0$, it suffices to show that the Albanese morphism $\alpha:X\to\Alb(X)$ factors through $f$.
By the universal property we know that$f$ factors through $\alpha$, say $f=g\alpha$ for some  surjective homomorphism $g:\Alb(X)\to A$ over $k$.
Fix, now, a generic point $b\in A$ over $k$, and consider $\alpha_b:f^{-1}(b)\to g^{-1}(b)$.
Our assumption on $f^{-1}(b)=:X_b$ is that there is, over some algebraically closed field $F$ extending $k(b)$, a birationally equivalent variety which is an algebraic homogeneous space for some connected linear algebraic group $L$ over $F$.
In particular, there exists a dominant rational map $L\to X_b$ over $F$, and hence $L\to \alpha_b(X_b)\subseteq A$.
But abelian varieties have no nontrivial rational subvarieties, and so we must have that $\alpha_b(X_b)$ is a point.
That is, $g:\Alb(X)\to A$ has a rational, and hence regular, section $h:A\to\Alb(X)$ such that $\alpha=hf$, as desired.

Suppose now that $v:X\to TX$ is an algebraic vector field such that $(X,v)$ admits no nontrivial rational first integrals.
Buium's theorem tells us that $(X,v)$ is isotrivial.
Working in a saturated model $\U$ of $\dcf$ with field of constants $\C$ extending $k$, and taking $p\in S(k)$ to be the generic type of $(X,v)$ over~$k$, we have that $p$ is $\C$-internal and weakly $\C$-orthogonal.
We have an abelian reduction $p\to p_{\ab}$, by Theorem~\ref{abredexists}.

Note that $p_{\ab}$ is $\C$-internal and weakly $\C$-orthogonal type over an algebraically closed field $k$ of constant parameters, with a commutative binding group.
So, by Fact~\ref{alldlog}, we may take $p_{\ab}$ to be the generic type of a logarithmic-differential equation on a connected abelian variety $A$ over $k$.
In particular, $p$ is the generic type of a $D$-variety $(A,w)$ over $k$.
The abelian reduction $p\to p_{\ab}$ gives rise to a dominant rational map of algebraic vector fields $f:(X,v)\to(A,w)$.
Since $X$ is a smooth projective variety and $A$ is a connected abelian variety, $f:X\to A$ is a surjective morphism.

Fix $a\models p$, $b:=f(a)\models p_{\ab}$, and $r:=\tp(a/kb)$.
Then $X_b$ is an irreducible $D$-subvariety of $(X,v)$ over $k(b)$ and $r$ is its generic type.
The fact that $r$ is stationary implies that $X_b$ is absolutely irreducible.

It remains to show that $X_b$ is birationally equivalent to an algebraic homogeneous space for a connected linear algebraic group.
Note that $r$ has trivial abelian reduction by Proposition~\ref{abfibab}, and so $\aut(r/\C)$, or rather its connected component,  is equal to its linear part by Proposition~\ref{trivabred}.
Since $\aut(r/\C)$ is a definable group over $k(b)$,  the main result of~\cite{kowalski2006quantifier} says that we may write 
$\aut(r/\C) = (G,v_G)^\delta$ where $(G,v_G)$ is a $D$-group over the (nonconstant) differential field $k(b)$.
The fact that $\aut(r/\C)$ is isomorphic to the $\mathcal C$-points of a linear algebraic group means that there exist a differential field extension $F\supseteq k(b)$ and $H$ a linear algebraic group over the constant field of $F$ such that $(G,v_G)_F \simeq (H,0)_F$.
In particular,  $G_F \simeq H_F$ is a linear algebraic group and hence so is $G$.
Now, the generically defined and generically transitive action $(G,v_G)^\delta \times (X_b,v)^\delta \dashrightarrow (X_b,v)^\delta$ can be lifted into  a generically defined and generically transitive action 
$(G,v_G) \times (X_b,v) \dashrightarrow (X_b,v)$
in the category of $D$-varieties over $k(b)$, see for example~\cite[Proposition~1.2]{HrIt}.
Forgetting about the differential structure, we obtain a generically defined generically  transitive action $G \times X_b \dashrightarrow X_b$ in the category of varieties over $k(b)$.
It then follows by Weil's group chunk theorem that $X_b$ is birationally equivalent to an algebraic homogeneous space under the action of $G$. 
\end{proof}

\bigskip
\section{Applications to compact complex manifolds}
\label{sec:ccm}

\noindent
There is a well known analogy between $\dcf$ and the theory of compact complex manifolds.
We should expect our work on abelian reductions to have bimeromorphic consequences analogous to the differential-algebraic ones discussed in the previous section.
We fulfill that expectation now.

Compact complex manifolds are studied model-theoretically by working in the theory $T=\ccm$ of the multistorted structure where there is a sort for each reduced and irreducible compact complex analytic space and a predicate for each complex analytic subset of each finite cartesian product of sorts.
The theory has many nice properties: quantifier elimination, elimination of imaginaries, and finiteness of Morley rank.
Moreover, a pure algebraically closed field is $0$-definable in this theory, which in the standard model is the complex field living on the projective line, and every nonmodular minimal type is nonorthogonal to it.
We do not review here how the main notions of geometric stability theory play out in $\ccm$; this material is covered in various accounts in the literature, to which we will make specific references as needed.
For now, let us point to the surveys~\cite{ccs-survey} and~\cite{moosapillay-survey}.

We work in a sufficiently saturated model $\U\models\ccm$ with $\PP$ the sort of the projective line and $\C$ the algebraically closed field definable on it.
See~\cite[$\S$2]{ret} for a discussion of the nonstandard model, in particular the description of  types in $\U$ as the generic types of nonstandard analytic sets.

For complex geometry we will follow the terminology of Fujiki~\cite{fujiki83}, which we now review.
By a {\em complex variety} we mean a reduced and irreducible complex analytic space.
A compact complex variety is {\em Moishezon} if it is the holomorphic (equivalently meromorphic) image of a projective algebraic variety.
A compact complex variety is of {\em K\"ahler-type} if it is the holomorphic (equivalently meromorphic) image of a compact K\"ahler manifold.

Suppose $X$ and~$Y$ are compact complex varieties and $f:X\to Y$ is a surjective morphism.
We say that a property {\em holds of the general fibre of $f$} if there is a proper analytic subset $E\subseteq Y$ such that the property holds of $X_y$ for all $y\in Y\setminus E$.
If, instead, we allow $E$ to be a countable union of proper analytic subsets, then this is what it means for the property to {\em hold of the `general' fibre}.
The morphism $f$ is called a {\em fibre space} if the general fibre is irreducible.
Every stationary type in $\ccm$ over a finite set of parameters is the generic type of a generic fibre of a fibre space.

Finally, let us recall base extensions of fibre spaces, which corresponds to taking nonforking extensions of stationary types.
Given a fibre space $f:X\to Y$ and a surjective morphism $g:Z\to Y$ where $Z$ is another compact complex variety, there is a unique maximal dimensional irreducible component of the fibred product $X\times_YZ$ that projects onto $Z$, which we denote by $X_{(Z)}$.
Moreover, the projection $f_{(Z)}:X_{(Z)}\to Z$ is itself a fibre space with general fibre the same as that of $f$.
See, for example,~\cite[Remark~2.4]{ret}, for a proof of these facts.
We call $f_{(Z)}$ the {\em base extension} of $f$ by $g$.
The model-theoretic content here is that if $b\in Z(\U)$ is a generic point and $p$ is the generic type of the fibre of $f$ over $g(b)$ then the generic type of the fibre of $f_{(Z)}$ over $b$ is the nonforking extension of $p$.
When $f=g$ we write $X_Y^{(2)}$ for the base extension, rather than $X_{(X)}$.

\medskip
\subsection{Binding groups of autonomous types and minimality}
In our applications to differential-algebraic geometry we mostly restricted our attention to types~$p$ over constant parameters.
The main reason for this was that in the autonomous case we have Fact~\ref{linearbg} which tells us that when $p$ is $\C$-internal weakly $\C$-orthogonal, if $\aut(p/\C)$ is linear then it is commutative.
In $\ccm$, at least if we restrict to K\"ahler-type manifolds, it turns out that something even stronger is true.
The following is a model-theoretic consequence of work of Fujiki on {\em relative Albanese} maps:

\begin{proposition}
\label{linearbg-ccm}
Suppose $B\subseteq\C$ and $p\in S(B)$ is a $\C$-internal weakly $\C$-orthogonal type living in the sort of a compact complex variety of K\"ahler-type.
If $p$ is nonalgebraic then $\aut(p/\C)$ is not linear.
\end{proposition}

\begin{proof}
Replacing $p$ by its unique extension $\overline p$ to $\acl(B)$, we may assume that $B=\acl(B)$.
By Theorems~\ref{abredexists} and~\ref{abred=alb} we have an abelian reduction $p\mapsto p_{\ab}$ with $\aut(p_{\ab}/\C)$ being the abelian part of $\aut(p/\C)$.
It therefore suffices to show that if~$p$ is nonalgebraic then $p$ has nontrivial abelian reduction.

The canonical base of $p$ is the definable closure of some finite tuple $b$ from~$B$.
So $q:=p\upharpoonright_b$ is $\C$-internal and weakly $\C$-orthogonal, and $p$ is the nonforking extension of $q$ to $B$.
It follows by Proposition~\ref{abred-nf} that $p_{\ab}$ is the nonforking extension of $q_{\ab}$ to $B$.
It therefore suffices to show that if $q$ is nonalgebraic then $q$ has nontrivial abelian reduction.

Since $q$ is a stationary type over finitely many parameters in $\C$, it is the generic type in $\U$ of the generic fibre of a fibre space $f:X\to Y$ where $Y$ is Moishezon.
As~$p$, and hence $q$, live on the sort of a compact complex variety of K\"ahler-type, we can take $X$ to be such.
The fact that $q$ is $\C$-internal implies that the general fibre of $f$ is Moishezon (see~\cite[Proposition~4.4]{ret}).
Under these conditions -- namely, that $X$ is of K\"ahler-type, $f$ is a fibre space, and the general fibres of $f$ are Moishezon -- Fujiki proves in~\cite[Theorem~2]{fujiki-alb} that a relative Albanese of $f$ exists; namely a fibre space $\Alb(X/Y)\to Y$ whose general fibre is an abelian variety, and a meromorphic map $\alpha:X\to\Alb(X/Y)$ over $Y$ such that, for general $y\in Y$, $\alpha_y:X_y\to\Alb(X_y)$ is the Albanese morphism of $X_y$.
Moreover, he shows that $\alpha$ is itself {\em Moishezon} in the sense that it embeds meromorphically into a projective linear space over $\Alb(X/Y)$.

Suppose now that  $q$ has trivial abelian reduction.
Then the connected component of the binding group, $\aut^\circ(q/\C)$, is equal to its linear part by Proposition~\ref{trivabred}.
The generic fibre of $f$ in $\U$ thus admits the action of a linear algebraic group with Zariski dense orbit.
The same is therefore true of the general fibre of $f$ in the standard model.
But this means that the Albanese of the general fibre of $f$ is trivial.
Hence $\Alb(X/Y)\to Y$ is a bimeromorphism.
That is, $\alpha$ is bimeromorphic to $f$ itself.
But then, as $\alpha$ is Moishezon, it follows that $f$ is Moishezon.
Since $Y$ is Moishezon, this in turn implies that $X$ is Moishezon.
The only way that can be consistent with $q$ being weakly $\C$-orthogonal is if $q$ is an algebraic type.
\end{proof}

This analogue of Fact~\ref{linearbg} allows us to extend to $\ccm$ two recent results on autonomous types in $\dcf$, that we take this opportunity to record:

\begin{corollary}
\label{autccm}
Suppose $B\subseteq\C$ and $p\in S(B)$ lives on the sort of a compact complex variety of K\"ahler-type.
Then the following hold:
\begin{itemize}
\item[(a)]
The {\em degree of nonminimality} of $p$ is at most $1$.
That is, either $p$ is minimal or for $a\models p$ there is a nonalgebraic forking extension of $p$ over $Ba$.
\item[(b)]
If every pair of distinct realisations of $p$ is independent over $B$ then $p$ is minimal.
\end{itemize}
\end{corollary}

\begin{proof}
It is shown in~\cite{nmdeg<2} that if $T$ is a complete theory satisfying:
\begin{enumerate}
\item
$T$ is a totally transcendental theory eliminating imaginaries and the ``there exists infinitely many" quantifier,
\item
there is a pure $0$-definable algebraically closed field $F$ to which every non locally modular minimal type is nonorthogonal, and
\item
if $q$ is a nonalgebraic type over parameters in $F$, and $q$ is internal and weakly orthogonal to $F$, then the binding group of $q$ has nontrivial center,
\end{enumerate}
then all types of finite rank in $T$ have degree of nonminimality at most~$1$.
So to prove part~(a), we take as $T$ the reduct of $\ccm$ to those sorts that are of K\"ahler-type.
Conditions~(1) and~(2) are then satisfied with $F=\C$.
It remains to verify~(3).
Suppose $B\subseteq\C$ and $q\in S(B)$ is nonalgebraic, $\C$-internal, and weakly $\C$-orthogonal.
Now $\aut(q/\C)$ is definably isomorphic to the $\C$-points of an algebraic group, and by~\cite[Theorem 13]{Ros} centerless algebraic groups are linear.
But Proposition~\ref{linearbg-ccm} ensures that $\aut(q/\C)$ is not linear.
Hence the binding group of $q$ has nontrivial center, as desired.

Part~(b) follows directly from part~(a) as follows:
If $p$ is not minimal then by part~(a) for every $a\models p$ we have a nonalgebraic forking extension $q$ of $p$ over $S(Ba)$.
Let $b\models q$.
By nonalgebraicity of $q$, $a$ and $b$ are distinct realisations of $p$. But they are not independent as $q$ is a forking extension.
\end{proof}

\begin{remark}
Part~(b) of Corollary~\ref{autccm} is the analogue in $\ccm$ of a result in $\dcf$ that appears in Theorem~C of~\cite{c3c2}.
This was established first.
It was only later, in the main theorem of~\cite{nmdeg<2}, that the $\dcf$ precursor of part~(a) was obtained.
\end{remark}

The above Corollary has a formulation in bimeromorphic geometry that we think is of independent interest, and that we now discuss.
Recall that in bimeromorphic geometry a compact complex variety $X$ is called {\em simple} if it cannot be covered by a family of proper infinite analytic subsets that cover $X$.
Here is a natural relativisation of this notion:

\begin{definition}
Suppose $f:X\to Y$ is a fibre space.
We will say that $f$ is {\em simple} if, for any surjective morphism $g:Z\to Y$, any irreducible proper analytic subset $W\subseteq X_{(Z)}$ that projects onto both $X$ and $Z$ must be of dimension $\dim Z$.
\end{definition}

Let us explain this a bit more.
Suppose we are given $W\subsetneq X_{(Z)}$ that projects dominantly onto both $X$ and $Z$ but with $\dim W>\dim Z$.
For general $y\in Y$, consider the fibre $W_y\subseteq X_y\times Z_y$.
This is a family of subsets of $X_y$ parametrised by $Z_y$ that covers $X_y$ (as $W$ projects onto $X$), and whose general members are proper (as $W\neq X_{(Z)}$), infinite (as $\dim W>\dim Z$), analytic subsets of $X_y$.
So $f:X\to Y$ being simple rules out the existence of covering families of proper infinite analytic subsets of $X_y$ given {\em uniformly} in $y\in Y$.
This is, in general, a strictly weaker condition than asking that the general (or even `general') fibre of~$f$ be simple.
Indeed, an old example of Lieberman~\cite{lieberman} studied by Campana~\cite{campana-exemples} and then model-theoretically by Pillay and Scanlon~\cite{pillayscanlon} gives an example where the general fibres are not simple but the fibration is in the above sense, see~\cite[$\S$2]{pillayscanlon}.
On the other hand, using the techniques of~\cite{saturated} on {\em essential saturation}, one can show that if $X$ is of K\"ahler-type (Lieberman's example is not) then $f$ being simple is equivalent to the `general' fibre being simple.

Simplicity of a fibre space $f$ is equivalent to minimality of the generic type of the generic fibre of $f$.

The following criterion for simplicity is a direct translation of Corollary~\ref{autccm} that we leave the reader to perform.
It says that as long as $X$ is of K\"ahler-type and $Y$ is Moishezon, to verify simplicity one need only consider the case of $f=g$ itself.

\begin{theorem}
\label{autccm-geometric}
Suppose $X$ is a compact complex variety of K\"ahler-type, $Y$ is Moishezon, and $f:X\to Y$ is a fibre space.
If every proper irreducible analytic subset of $X_Y^{(2)}$ that projects onto $X$ in both co-ordinates is of dimension $\dim X$ then $f$ is simple.
\end{theorem} 
 
\medskip
\subsection{Bounding nonorthogonality in $\ccm$}
The following is the bimeromorphic analogue of a vector field admitting no nontrivial rational first integrals:

\begin{definition}
Suppose $f:X\to Y$ is a fibre space.
Let us say that $f$ {\em admits no nonconstant meromorphic functions} if any meromorphic function on $X$ is constant on the general fibre of $f$.
\end{definition}

\begin{remark}
Note that this is equivalent to saying that $f$ induces an isomorphism of the meromorphic function fields of $X$ and $Y$.
On the other hand it is strictly weaker than asking that the {\em algebraic dimension of~$f$} be zero.
The latter would mean that the `general' fibre of $f$ admits no nonconstant meromorphic functions.
But we are only ruling out {\em uniformly} given nonconstant meromorphic functions on the fibres.
For example, if $f:X\to Y$ is an {\em algebraic reduction} of $X$, namely $Y$ is Moishezon and $f$ induces an isomorphism between meromorphic function fields, then $f$ admits no nonconstant meromorphic functions.
However, the fibres of $f$ may very well themselves be Moishezon and positive-dimensional, having therefore many nonconstant meromorphic functions.
\end{remark}

Given a fibre space $f:X\to Y$ 
We denote by $f^{(2)}:X_Y^{(2)}\to Y$ the fibre space obtained by composing $X_Y^{(2)}\to X$ with $f$ itself.
So the general fibre of $f^{(2)}$ is the cartesian square of the general fibre of $f$.
Similarly, one can define $f^{(n)}$ for all $n\geq 1$.
If $p$ is the generic type of the generic fibre of $f$ then $p^{(n)}$ is the generic type of the generic fibre of $f^{(n)}$.

The following application of abelian reductions to bimeromorphic geometry is the analogue of Corollary~\ref{nwodeg-geometric} above.

\begin{theorem}
\label{nwodeg-ccm}
Suppose $X$ is a compact complex variety of K\"ahler-type, $Y$ is Moishezon, and $f:X\to Y$ is a fibre space.
The following are equivalent:
\begin{itemize}
\item[(i)]
 $f^{(2)}$ admits no nonconstant meromorphic functions,
 \item[(ii)]
 $f^{(n)}$ admits no nonconstant meromorphic functions, for all~$n\geq 1$, and
 \item[(iii)]
 no base extension of $f$ admits a nonconstant meromorphic function.
 \end{itemize}
\end{theorem}

\begin{proof}
Let $a\in X(\U)$ be a generic point of $X$, so that $b:=f(a)\in Y(\U)$ is generic in~$Y$.
Let $p:=\tp(a/b)$.
The fact that $f$ is a fibre space implies that $p$ is stationary, see~\cite[Lemma~2.7]{ret}.
That $Y$ is Moishezon means we can take $b$ to be a finite tuple from $\C$.
The property of the fibre space $f$ admitting no nonconstant meromorphic functions expresses precisely that the type $p$ is weakly $\C$-orthogonal.
Hence the equivalence of parts~(ii) and~(iii) is just the equivalence of $p^{(n)}$ being weakly $\C$-orthogonal for all $n\geq 1$ with $p$ being $\C$-orthogonal.
And of course, that~(ii) implies~(i) is clear.
It will suffice, therefore, to show that~(i) implies~(iii).
That is, if $p$ is not $\C$-orthogonal then $p^{(2)}$ is not weakly $\C$-orthogonal.
But this is precisely Theorem~\ref{nwodeg} in the case of $\dcf$.
Inspecting that proof we see that everything goes through except for the final appeal to Fact~\ref{linearbg}, which we replace by Proposition~\ref{linearbg-ccm}.
\end{proof}

One can ask about the $\ccm$ analogues of the other applications to algebraic vector fields given in Section~\ref{sec:dcf}.
But these lead to results that are already known, and that were, indeed, the original inspiration for looking at abelian reductions model-theoretically.
So, for example, the classification of $2$-dimensional isotrivial algebraic vector fields admitting no nontrivial rational first integrals, given in Corollary~\ref{classification-ordertwo} above, corresponds in the case of bimeromorphic geometry to Fujiki's classification in~\cite{fujiki83} of $3$-dimensional compact K\"ahler manifolds of algebraic dimension~$1$.
The analogue of Theorem~\ref{nofibre-dcf-geometric} above, on the structure types over constant parameters that admit no proper almost-fibrations, also exists in $\ccm$, but it's geometric formulation appears already as Propositions~4.2 and~4.3 of~\cite{moosa2014some}, but see also the earlier~\cite[Theorem~2.4]{COP}.
Finally, Theorem~\ref{albanese-surjective} above, which describes the Albanese map of a smooth projective variety admitting an algebraic vector field with no nontrivial rational first integrals, corresponds to Theorem~1 of~\cite{fujiki83}.

\bigskip

\begin{thebibliography}{10}

\bibitem{brion}
Michel Brion.
\newblock Some structure theorems for algebraic groups.
\newblock In {\em Algebraic groups: structure and actions}, volume~94 of {\em
  Proc. Sympos. Pure Math.}, pages 53--126. Amer. Math. Soc., Providence, RI,
  2017.

\bibitem{Buium}
Alexandru Buium.
\newblock {\em Differential function fields and moduli of algebraic varieties}.
\newblock Springer-Verlag, 1986.

\bibitem{campana-exemples}
Fr\'{e}d\'{e}ric Campana.
\newblock Exemples de sous-espaces maximaux isol\'{e}s de codimension deux d'un
  espace analytique compact.
\newblock In {\em Institut \'{E}lie {C}artan, 6}, volume~6 of {\em Inst.
  \'{E}lie Cartan}, pages 106--127. Univ. Nancy, Nancy, 1982.

\bibitem{COP}
Fr\'{e}d\'{e}ric Campana, Keiji Oguiso, and Thomas Peternell.
\newblock Non-algebraic hyperk\"{a}hler manifolds.
\newblock {\em J. Differential Geom.}, 85(3):397--424, 2010.

\bibitem{c3c2}
James Freitag, R\'emi Jaoui, and Rahim Moosa.
\newblock When any three solutions are independent.
\newblock {\em Invent. Math.}, 2022, no. 3, p.~1249-1265.

\bibitem{nmdeg<2}
James Freitag, R\'{e}mi Jaoui, and Rahim Moosa.
\newblock The degree of nonminimality is at most 2.
\newblock {\em J. Math. Log.}, 23(3):Paper No. 2250031, 6, 2023.

\bibitem{nmdeg}
James Freitag and Rahim Moosa.
\newblock Bounding nonminimality and a conjecture of {B}orovik-{C}herlin.
\newblock {\em Journal of the European Mathematical Society}, 2023.
\newblock Online first, DOI 10.4171/JEMS/1384.

\bibitem{fujiki83}
Akira Fujiki.
\newblock On the structure of compact complex manifolds in {C}.
\newblock In {\em Algebraic varieties and analytic varieties ({T}okyo, 1981)},
  volume~1 of {\em Adv. Stud. Pure Math.}, pages 231--302. North-Holland,
  Amsterdam, 1983.

\bibitem{fujiki-alb}
Akira Fujiki.
\newblock Relative algebraic reduction and relative {A}lbanese map for a fiber
  space in {${\mathcal C}$}.
\newblock {\em Publ. Res. Inst. Math. Sci.}, 19(1):207--236, 1983.

\bibitem{furstenberg}
H.~Furstenberg.
\newblock {\em Recurrence in ergodic theory and combinatorial number theory}.
\newblock Princeton University Press, Princeton, N.J., 1981.
\newblock M. B. Porter Lectures.

\bibitem{hrushovski1989almost}
Ehud Hrushovski.
\newblock Almost orthogonal regular types.
\newblock {\em Annals of Pure and Applied Logic}, 45(2):139--155, 1989.

\bibitem{udicomputing}
Ehud Hrushovski.
\newblock Computing the {G}alois group of a linear differential equation.
\newblock In {\em Differential {G}alois theory ({B}edlewo, 2001)}, volume~58 of
  {\em Banach Center Publ.}, pages 97--138. Polish Acad. Sci. Inst. Math.,
  Warsaw, 2002.

\bibitem{HrIt}
Ehud Hrushovski and Masanori Itai.
\newblock On model complete differential fields.
\newblock {\em Transactions of the American Mathematical Society}, Volume
  355(11):4267--4296, 2003.

\bibitem{HrSo}
Ehud Hrushovski and \v{Z}eljko Sokolovi\'{c}.
\newblock Strongly minimal sets in differentially closed fields.
\newblock {\em unpublished manuscript}, 1993.

\bibitem{jaoui20}
R\'{e}mi Jaoui.
\newblock Corps diff\'{e}rentiels et flots g\'{e}od\'{e}siques
  {I}---{O}rthogonalit\'{e} aux constantes pour les \'{e}quations
  diff\'{e}rentielles autonomes.
\newblock {\em Bull. Soc. Math. France}, 148(3):529--595, 2020.

\bibitem{JJP}
R\'{e}mi Jaoui, L\'{e}o Jimenez, and Anand Pillay.
\newblock Relative internality and definable fibrations.
\newblock {\em Adv. Math.}, 415:Paper No. 108870, 38, 2023.

\bibitem{leo}
L\'{e}o Jimenez.
\newblock Groupoids and relative internality.
\newblock {\em J. Symb. Log.}, 84(3):987--1006, 2019.

\bibitem{jin-moosa}
Ruizhang Jin and Rahim Moosa.
\newblock Internality of logarithmic-differential pullbacks.
\newblock {\em Trans. Amer. Math. Soc.}, 373(7):4863--4887, 2020.

\bibitem{KolchinDAAG}
Ellis~R. Kolchin.
\newblock {\em Differential Algebra and Algebraic Groups}.
\newblock Academic Press, New York, 1976.

\bibitem{kowalski2006quantifier}
Piotr Kowalski and Anand Pillay.
\newblock Quantifier elimination for algebraic {$D$}-groups.
\newblock {\em Trans. Amer. Math. Soc.}, 358(1):167--181, 2006.

\bibitem{lieberman}
David~I. Lieberman.
\newblock Compactness of the {C}how scheme: applications to automorphisms and
  deformations of {K}\"{a}hler manifolds.
\newblock In {\em Fonctions de plusieurs variables complexes, {III} ({S}\'{e}m.
  {F}ran\c{c}ois {N}orguet, 1975--1977)}, volume 670 of {\em Lecture Notes in
  Math.}, pages 140--186. Springer, Berlin, 1978.

\bibitem{ret}
Rahim Moosa.
\newblock A nonstandard {R}iemann existence theorem.
\newblock {\em Trans. Amer. Math. Soc.}, 356(5):1781--1797, 2004.

\bibitem{ccs-survey}
Rahim Moosa.
\newblock The model theory of compact complex spaces.
\newblock In {\em Logic {C}olloquium '01}, volume~20 of {\em Lect. Notes Log.},
  pages 317--349. Assoc. Symbol. Logic, Urbana, IL, 2005.

\bibitem{saturated}
Rahim Moosa.
\newblock On saturation and the model theory of compact {K}\"{a}hler manifolds.
\newblock {\em J. Reine Angew. Math.}, 586:1--20, 2005.

\bibitem{moosapillay-survey}
Rahim Moosa and Anand Pillay.
\newblock Model theory and {K}\"{a}hler geometry.
\newblock In {\em Model theory with applications to algebra and analysis.
  {V}ol. 1}, volume 349 of {\em London Math. Soc. Lecture Note Ser.}, pages
  167--195. Cambridge Univ. Press, Cambridge, 2008.

\bibitem{moosa2014some}
Rahim Moosa and Anand Pillay.
\newblock Some model theory of fibrations and algebraic reductions.
\newblock {\em Selecta Mathematica}, 20(4):1067--1082, 2014.

\bibitem{GST}
Anand Pillay.
\newblock {\em Geometric Stability Theory}.
\newblock Oxford University Press, 1996.

\bibitem{pillay2004}
Anand Pillay.
\newblock Algebraic {$D$}-groups and differential {G}alois theory.
\newblock {\em Pacific J. Math.}, 216(2):343--360, 2004.

\bibitem{pillayscanlon}
Anand Pillay and Thomas Scanlon.
\newblock Compact complex manifolds with the {DOP} and other properties.
\newblock {\em J. Symbolic Logic}, 67(2):737--743, 2002.

\bibitem{Ros}
Maxwell Rosenlicht.
\newblock Some basic theorems on algebraic groups.
\newblock {\em Amer. J. Math.}, 78:401--443, 1956.

\bibitem{umemura}
Hiroshi Umemura.
\newblock Second proof of the irreducibility of the first differential equation
  of {P}ainlev\'{e}.
\newblock {\em Nagoya Math. J.}, 117:125--171, 1990.

\end{thebibliography}

\end{document}